\documentclass[leqno,CJK]{siamltex704}
\usepackage{amssymb,amsmath,graphicx,amscd,mathrsfs}
\usepackage{color,xcolor,amsmath}
\usepackage{amsmath}
\usepackage{graphicx}
\usepackage{epstopdf}
\usepackage{subfigure}
\usepackage{mathrsfs}
\usepackage{float}
\usepackage{amsfonts,amssymb}
\usepackage{dsfont}
\usepackage{pifont}
\usepackage{hyperref}
\usepackage{multirow}
\usepackage{placeins}
\numberwithin{equation}{section}
\def\3bar{{|\hspace{-.02in}|\hspace{-.02in}|}}

\def\b0{\boldsymbol{0}}

\def\bf{{\mathbf{f}}}

\newtheorem{defi}{Definition}[section]

\newtheorem{assumption}{Assumption}[section]
\newtheorem{algorithm1}{Weak Galerkin Algorithm}

 \newcommand{\Real}{\mathbb{R}}

\allowdisplaybreaks 

\begin{document}

\title{Weak Galerkin methods for the Stokes eigenvalue problem}

\author{
Yunying Fan
\and
Qilong Zhai
}
\maketitle
\begin{abstract}
In this paper, we rewrite the Stokes eigenvalue problem as an Elliptic eigenvalue problem restricted to subspace, and introduce an abstract framework of solving abstract elliptic eigenvalue problem from \cite{zhai2019weak} to give the WG scheme, error estimates and asymptotic lower bounds. Besides, we introduce a new stabilizer and several inequalities from \cite{carstensen2020skeletal} to prove GLB properties. Some numerical examples are provided to validate our theoretical analysis.
\end{abstract}



\section{Introduction}

The eigenvalue problem arises from many {\color{blue}applications} of \\mathematics and
physics, including quantum mechanics, fluid mechanics, stochastic process
and structural mechanics. Eigenvalue problem, especially, the Stokes eigenvalue problem, has attracted much interest since it’s frequently encountered in a variety of applications, for instance, to study the stability of N-S equation, and also appears in the analysis of the elastic stability of thin plates \cite{osborn1976approximation,mercier1981eigenvalue} and so on. Many numerical methods have been developed for the Stokes eigenvalue problem, such as finite difference methods \cite{da2009mimetic}, finite element methods \cite{dohrmann2004stabilized,gedicke2018arnold} and finite volume method \cite{eymard2006stabilized}.

The finite element method is one of the efficient approaches for the Stokes eigenvalue problem for its simplicity and adaptivity on triangular meshes. At the same time, since the all eigenvalues of Stokes equation is real, we can get accurate intervals of eigenvalues by the upper bounds and lower bounds. However, due to the minimum-maximum principle, the conforming finite element method always gives the upper bounds \cite{1965Difference,strang1974analysis}.  In order to get accurate intervals for eigenvalues, it is necessary to have lower bounds of eigenvalues. There are mainly two ways, the post-processing method \cite{Larson2000APA,Liu2013VerifiedEE} and the nonconforming finite element method \cite{Lin2008NewEO}. The main difficulty of the first way is an auxiliary problem must be solved and may suffer reduced convergence order. The second way can provide asymptotic lower bounds for eigenvalues without solving an auxiliary problem, while it seems difficult to construct a high order or high dimension nonconforming element.

But recently, a new method for solving the PDE, named the weak Galerkin (WG) method, can solve those problem.  The main features of WG method are that replace differential operator with weak differential operator and choose discontinuous piecewise polynomials as base function. Hence, it can be applied to many types of areas and partitions, and is very flexible and easy to expand. The WG method was first introduced in \cite{Wang2014a} for the second order elliptic equation, and was soon applied to many types of partial differential equations, such as the Stokes equation \cite{Wang2016AWG}, the parabolic equation \cite{Li2013}, the biharmonic equation\cite{Oden2007,Mu2013,Zhang2015}, the Brinkman equation \cite{mu2014stable}, and the Maxwell equation \cite{Mu2013d}. In \cite{Xie2015}, the Laplacian eigenvalue problem was investigated by the WG method. And in \cite{zhai2019weak}, a theoretical framework was designed to solve the general abstract elliptic eigenvalue problem by the WG method. It is this framework that we used to investigate by the WG method. An {\color{blue}interesting} feature is: it offers asymptotic lower bounds for the Laplacian eigenvalues or general abstract elliptic eigenvalues problem on polygonal meshes by employing high order polynomial elements.

However, asymptotic lower bound is not what we want because it also suffers reduced convergence order. In \cite{carstensen2020skeletal}, the Laplacian eigenvalue problem was given the conditions of GLB by skeletal finite element method. Since the similarity between skeletal eigenvalue method and WG method, we can give a rough proof of GLB properties of Stokes eigenvalue problem by similar method used in \cite{carstensen2020skeletal}.

An outline of the paper goes as follows. In Section 2, we introduce the equivalent form of Stokes eigenvalue problem and some notations. In Section 3, we introduce a theoretical framework of solving abstract elliptic eigenvalue problem and use it to give the WG scheme, error estimates and asymptotic lower bounds of Stokes eigenvalue problem in Section 4, section 5 and section 6, respectively. In Section 7, we introduce a new stabilizer and several inequalities to prove GLB properties of Stokes eigenvalue problem.In Section 8, some numerical examples are presented to {\color{blue}verify} our theoretical analysis. Some concluding remarks are given in the final section.


\section{Notations and Weak form of Stokes eigenvalue problem}
In this section, we state some notation and
introduce the standard WG scheme for Stokes eigenvalue 
Firstly,we introduce the notation as followings:$$W^{s,p}(\Omega)=\{u\in L^p(\Omega);\forall \kappa,|\kappa|\leq s,\text{weak derivatives }\partial^{\kappa}u\text{ exists, }\partial^{\kappa}u \in L^p(\Omega)\}$$
its associated norm:
$$\parallel u\parallel_{s,p}=(\sum_{|\kappa|\leq s}\parallel\partial^{\kappa}u\parallel_{L^p}^p)^{\frac{1}{p}}$$
For $p=2$, $W^{s,2}(\Omega)=H^s (\Omega)$ is Hilberts space. For $s=0$, $L^p(\Omega)=W^{0,p}(\Omega)$.
Especially,we introduce several following Hilbert space,
$$V=[H_0^1(\Omega)]^d(d=2,3),\quad Q=L_0^2(\Omega)=\{q\in L^2(\Omega):\int_{\Omega}q dx=0\}.$$

In this paper,we consider the following Stokes eigenvalue problem:
Find$(\lambda;\mathbf{u},p)  \in \mathbb{R}\times V \times Q$ such that
\begin{equation}
    \begin{cases}
        -\Delta \mathbf{u}+\nabla p=\lambda \mathbf{u},& \text{in}\Omega,\\
        \nabla \cdot \mathbf{u}=0,& \text{in}\Omega,\\ 
        \mathbf{u}=\mathbf{0},& \text{on}\partial\Omega, \\
        \int_{\Omega}^{} |\mathbf{u}|^2 d\Omega=1,
     \end{cases}
\end{equation}
where $\Omega \subset \mathbb{R}^d(d=2,3)$ denotesss the computing domain with $Lipschitz$ boundary $\partial\Omega$; $\mathbf{u}=\mathbf{u}(\mathbf{x})$ is the velocity $\mathbf{p}=\mathbf{p}(\mathbf{x})$ is pressure. In addition,symbols $\Delta , \nabla$and$\nabla \cdot$ denotesss the Laplacian, gradient and divergence operators,respectively.

For the aim of discretization,we give a weak form as follows:Find$ (\lambda,\mathbf{u},p)$   $\in \mathbb{R}\times V \times Q$ such that $b(\mathbf{u},\mathbf{u})=1$ and\\
\begin{equation}
    \begin{cases}
        a(\mathbf{u},\mathbf{v})-c(\mathbf{v},p)=\lambda b(\mathbf{u},\mathbf{v}), \forall \mathbf{v} \in V,\\
        c(\mathbf{u},q)=0,      \forall q \in Q,
     \end{cases}
\end{equation}\\
where$$a(\mathbf{u},\mathbf{v}):=\int_{\Omega}^{}\nabla \mathbf{u}:\nabla \mathbf{v} d\Omega,\quad c(\mathbf{v},p):=\int_{\Omega}^{}p\nabla \cdot \mathbf{v} d\Omega,$$
$$b(\mathbf{u},\mathbf{v}):=\int_{\Omega}^{}\mathbf{u}\cdot \mathbf{v} d\Omega,$$
and $\nabla \mathbf{u}:\nabla \mathbf{v}=\sum_{i=1}^{d}\sum_{j=1}^{d}\frac{\partial u_i}{\partial x_j}\frac{\partial v_i}{\partial x_j}.$

In this paper,we also will give an equivalent formulation of the above form (equivalent is visible in \cite{Xie2018ExplicitLB}), we define
$$V_0:=\{\mathbf{v}\in V:c(\mathbf{v},q)=0,\forall q\in Q\}.$$
 Then,we can get the formulation as following:Find$(\lambda,\mathbf{u})\in \mathbb{R}\times V_0$ such that $b(\mathbf{u},\mathbf{u})=1$ and
\begin{equation}
    a(\mathbf{u},\mathbf{v})=\lambda b(\mathbf{u},\mathbf{v}),\quad \forall \mathbf{v}\in V_0.
\end{equation}
From \cite{Xie2018ExplicitLB}, we know eigenvalue problems (2.2) and (2.3) have the same eigenvalue sequences $\{\lambda_i\}_{i=1}^{\infty}$
$$0<\lambda_1\leq\cdots\leq\lambda_i,\lim_{i\to\infty}\lambda_i=\infty,$$
and corresponding$$(\mathbf{u_1},p_1),\cdots,(\mathbf{u_i},p_i),\cdots,$$where $b(\mathbf{u_i},\mathbf{u_j})=\delta_{ij},$ and $\delta_{ij}$ denotesss Kronecker symbol.

Rewriting the Stokes eigenvalue problem into the variational form is convenient to solved by following framework of the weak finite element method for the abstract elliptic eigenvalue problem.

\section{Framework for abstract elliptic eigenvalue problems}
In this section, we introduce a general framework of WG  method from \cite{zhai2019weak} for abstract elliptic eigenvalue problems.

Consider the following abstract eigenvalue problem : Find$(\lambda,u)\in \mathbb{R}\times V_c$ such that $b(u,u)=1$ and
\begin{equation}
    a(u,v)=\lambda b(u,v),\quad \forall v\in V_c;
\end{equation}

And corresponding discrete problem:Find$(\lambda_h,u_h)\in \mathbb{R}\times V_h\text{ such that }b(u_h,u_h)=1$ and
\begin{equation}
    a_w (u_h ,v_h)=\lambda_h  b(u_h ,v_h),\quad \forall v_h\in V_h;
\end{equation}

where $V_c$ and $V_h$ are the subset of Hilbert space V, $\parallel\cdot\parallel_{V}$ denotesss norm introduced by inner product defined on Hilbert space V, $a(\cdot,\cdot)$, $a_w(\cdot,\cdot)$ and $b(\cdot,\cdot)$ denotesss bilinear forms defined on $V_c$, $V_h$ and V, respectively.

The abstract framework can accomplish the error estimation of the eigenvalues and eigenfunctions, and the corresponding asymptotic lower bound estimation, under the following seven assumptions. In this section, we present these seven assumptions, and the final error estimates of the eigenvalues and eigenfunctions, and the lower bound estimates, respectively. All conclusions and derivations in this section can be found in \cite{zhai2019weak}.
\begin{assumption}
 $a(\cdot,\cdot)$, $a_w(\cdot,\cdot)$ and $b(\cdot,\cdot)$ are symmetric defined on $V_c$,$V_h$ and V, satisfying \\
 $$a(v,v)\geq\gamma_c\parallel v\parallel_V^2,\quad \forall v \in V_c,$$
 $$a_w(v_h,v_h)\geq\gamma(h)\parallel v_h\parallel_V^2,\quad \forall v_h \in V_h,$$
 where $\gamma_c$is positive constant, $\gamma(h)$is positive function.
\end{assumption}

Eigenvalue problem can be viewed as operator spectrum problem. Define following two operators $K:V_c\rightarrow V_c$ and $K_h:V_h\rightarrow V_h$ satisfying:
$$a(Kf,v)=b(f,v),\quad \forall v\in V_c,$$
$$a_w(K_h f_h,v_h)=b(f_h,v_h),\quad \forall v_h\in V_h.$$

After we have assumption 3.1, from Lax-Milgram theorem, easy to know that operator K and $K_h$ are meaningful. Hence,we can give following assumption.
\begin{assumption}
operator $K$ and $K_h$ is compact.
\end{assumption}
\begin{assumption}
there exists a bounded and linear operator $Q_h:V\rightarrow V_h$ statisfying
$$Q_h v_h=v_h,\quad \forall v_h\in V_h,$$
$$b(Q_h w,v_h)=b(w,v_h),\quad \forall w\in V,\quad v_h\in V_h.$$
\end{assumption}
For eigenvalue $\mu$ of operator $K$ ,we can define following operator spectrum mapping:
$$E_{\mu}(K)=\frac{1}{2\pi i}\int_{\Gamma_{\mu}}R_z(K)dz.$$
Where $R_z(K)=(z I-K)^{-1}$ denotesss resolvent operator. $R(E_{\mu}(K))$denotesss the range of spectrum mapping $E_{\mu}(K)$ of operator $K$, and it also represents eigenfunction space according to $\mu$. Similarly , $R(E_{\mu,h}(K_h))$ denotesss the range of spectrum mapping $R(E_{\mu,h}(K_h))$ of operator $K_h$, and it also represents eigenfunction space according to $\mu$.

In eigenfunction space, we define following errors:
$$e_{h,\mu}=\parallel(K-K_h Q_h)|_{R(E_{\mu}(K))}\parallel_V,\quad \delta_{h,\mu}=\displaystyle \sup_{\mathop{}_{\parallel u\parallel_ V=1}^{u\in R(E_{\mu}(K))}} \parallel u-Q_h u\parallel_V$$
\begin{assumption}
when $h\rightarrow 0$, $e_{h,\mu}\rightarrow 0$, and $\delta_{h,\mu}\gamma(h)^{-1}\rightarrow 0.$
\end{assumption}

Let X and Y is subset of Banach space V, we define following distance:
$$\rho_V (x,Y)=\displaystyle\inf_{y\in Y}\parallel x-y\parallel_V,\quad \rho_V(X,Y)=\displaystyle\sup_{\mathop{}_{x\in X}^{\parallel x \parallel_V=1}} \rho_V(x,Y),$$
$$\rho_V^{\prime}(X,Y)=max\{\rho_V (X,Y),\rho_V(Y,X)\}$$
\begin{theorem}
Under the assumption we have,
$$\rho_V^{\prime}(R(E_{\mu}(K)),R(E_{\mu,h}(K_h)))\leq C(e_{h,\mu}+\delta_h\gamma(h)^{-1}).$$where C is positive function about $\mu$.
\end{theorem}

Let $X$ is a Hilbert space with inner product $b(\cdot,\cdot)$, $\parallel \cdot \parallel_X$ is corresponding norm. In elliptic eigenvalue problem, X is $L^2(\Omega)$. Let $\Pi_0$ is projection from V to X with inner product $b(\cdot,\cdot)$, hence we can define following bounded linear operator $\Pi_0 K|_{\Pi_0 V_c}:\Pi_0 V_c\rightarrow\Pi_0 V_c$ and $\Pi_0 K_h|_{\Pi_0 V_h}:\Pi_0 V_h\rightarrow\Pi_0 V_h$.
Hence,in eigenfunction space,we define following errors:
$$e_{h,\mu}^{\prime}=\parallel(\Pi_0 K-\Pi_0 K_h Q_h)|_{R(E_{\mu}(K))}\parallel_X,\quad \delta_{h,\mu}^{\prime}=\rho_X(V,V_h)$$
\begin{assumption}
When $h\rightarrow 0$,$e_{h,\mu}^{\prime}\rightarrow 0$, and $\delta_{h}^{\prime}\gamma(h)^{-1}\rightarrow 0.$
\end{assumption}
\begin{theorem}
Under the assumption, we have,
$$\rho_X^{\prime}(R(E_{\mu}(K)),R(E_{\mu,h}(K_h)))\leq C(e_{h,\mu}^{\prime}+\delta_h^{\prime}\gamma(h)^{-1}).$$where C is positive function about $\mu$.
\end{theorem}

Under the assumption,like lemma 2.9 in \cite{zhai2019weak}, we can get $\lambda-\lambda_h=a(u,u)-a_w(Q_h u,Q_h u)$
$+a_w(u_h-Q_h u,u_h-Q_h u)-\lambda_h b(u-u_h,u-u_h).$ To do further estimate, we need to calculate the first two items. Hence we take:
$\epsilon_{h,u}=a(u,u)-a_w(Q_h u,Q_h u)$,
And assume that,
\begin{assumption}
when $h\rightarrow 0$,$\epsilon_{h,u}\rightarrow 0$, for $\mu \in \sigma(K),u\in R(E_{\mu}(K))$.
\end{assumption}

Combining the above assumptions, we can derive the following eigenvalue estimates.
\begin{theorem}
Let $\lambda$ is an eigenvalue with multiplicity m, $\{\lambda_{h,j}\}_{j=1}^m$ is its m numerical solution. $\{u_j\}_{j=1}^m$ is a group eigenvalue vectors according to $\lambda$. Under the assumption, when h is enough small,
$$|\lambda-\lambda_{h,j}|\leq C(\epsilon_{h,u_j}+e_{h,\lambda^{-1}}^2+e_{h,\lambda^{-1}}^{\prime  2}+\delta_h^2\gamma(h)^{-2}+\delta_h^{\prime2}\gamma(h)^{-2}).$$
\end{theorem}
\begin{assumption} $(\lambda,u)$ and $(\lambda_h,u_h)$ are the exact solution and numerical solution of  (3.1) and
(3.2),then
$$\epsilon_{h,u}\geq \lambda_h \parallel u-u_h\parallel_{X}^2.$$
\end{assumption}
Above all, we have following result.
\begin{theorem}
$(\lambda,u)$ and $(\lambda_h,u_h)$ are the exact solution and numerical solution of  (3.1) and (3.2), then $$\lambda\geq\lambda_h.$$
\end{theorem}

\section{WG scheme}

In this section, we apply the theoretical framework in Section 3 to give the WG space and numerical scheme of WG method for the Stokes eigenvalue problem.
For the following Stokes eigenvalue problem:
\begin{equation}
    \begin{cases}
        -\Delta \mathbf{u}+\nabla p=\lambda \mathbf{u},& \text{in }\Omega,\\
        \nabla \cdot \mathbf{u}=0,& \text{in }\Omega,\\ 
        \mathbf{u}=\mathbf{0},& \text{on }\partial\Omega, \\
        \int_{\Omega}^{} |\mathbf{u}|^{2} d\Omega=1,
     \end{cases}
\end{equation}
where $\Omega$ is a polygon region in $\Real^d$ $(d=2,3)$. 
From the analysis of section two, we know its equivalent variational form:Find $(\lambda,\mathbf{u})\in \mathbb{R}\times V_0$ such that $b(\mathbf{u},\mathbf{u})=1$ and
\begin{equation}
    a(\mathbf{u},\mathbf{v})=\lambda b(\mathbf{u},\mathbf{v}),\quad \forall \mathbf{v}\in V_0.
\end{equation}
where $a(\mathbf{u},\mathbf{v})=(\nabla \mathbf{u},\nabla \mathbf{v}), b(\mathbf{u},\mathbf{v})=(\mathbf{u},\mathbf{v}), c(\mathbf{v},q)=(\nabla \cdot \mathbf{v},q), V_0:=\{\mathbf{v}\in [H_0^1 (\Omega)]^d:c(\mathbf{v},q)=0,\forall q\in L_0^2 (\Omega)\}.$

Since $V_0$ is Hilbert space, this form satisfy the requirement of section three. Hence, we will give the WG scheme of problem (4.2) ,then verify hypothesis by hypothesis to get the error estimates and asymptotic lower bound estimates. In addition, we always assume that C is constants independent of grid size, We generally use $a\lesssim b$ to denotes $a\leq C b$.

In the following, we introduce some notation in the WG scheme, let $T_h$  a mesh partition of the region $\Omega$, and each of the mesh cells is a polygon or polyhedron satisfying the regularity assumption(details to see \cite{mu2012weak}). We use $h=\displaystyle\max_{T\in T_h}h_T$ to denotes the mesh size, $\epsilon_h$ to denotes the middle edge of $T_h$ and $\epsilon_h^0$ to denotes the inner edge $\epsilon_h\backslash \partial \Omega$.

To define weak gradient, we need weak function $\mathbf{v}=\{\mathbf{v}_0,\mathbf{v}_b\}$ satisfying $\mathbf{v}_0 \in [L^2(T)]^d$ and $\mathbf{v}_b \in [L^2(\partial T)]^d$. Let $V(T)$ denotess weak vector-value function space, i.e.
$$V(T)=\{\mathbf{v}=\{\mathbf{v}_0,\mathbf{v}_b\}:\mathbf{v}_0\in [L^2(T)]^d,\mathbf{v}_b\in [L^2(\partial T]^d\}.$$
\begin{defi}
weak gradient operator is denoted by $\nabla_{\omega}$,and defined by a uniquely determined polynomial in  $[P_r(T)]^{d\times d}$,satisfying that
$$\left \langle\nabla_{\omega}\mathbf{v},\mathbf{q}\right \rangle_T=-(\mathbf{v}_0,\nabla \cdot \mathbf{q})_T+\left \langle \mathbf{v}_b,\mathbf{q}\cdot \mathbf{n}\right \rangle_{\partial T},\quad \forall \mathbf{q}\in [P_r(T)]^{d\times d}.$$
\end{defi}
To define weak divergence, we need weak function $\mathbf{v}=\{\mathbf{v}_0,\mathbf{v}_b\}$ satisfying $\mathbf{v}_0 \in [L^2(T)]^d$ and $\mathbf{v}_b\cdot \mathbf{n} \in [L^2(\partial T)]^d$. $V(T)$ denotes weak vector-value function space, 
i.e.
$$V(T)=\{\mathbf{v}=\{\mathbf{v}_0,\mathbf{v}_b\}:\mathbf{v}_0\in [L^2(T)]^d,\mathbf{v}_b \cdot \mathbf{n}\in [L^2(\partial T]^d\}.$$

\begin{defi}
Weak divergence operator is denoted by $\nabla_{\omega}\cdot $,and defined by a uniquely determined polynomial in $P_r(T)$,satisfying that
$$\left \langle\nabla_{\omega }\cdot \mathbf{v},\phi\right \rangle_T=-(\mathbf{v}_0,\nabla  \phi)_T+\left \langle \mathbf{v}_b\cdot \mathbf{n},\phi\right \rangle_{\partial T},\quad \forall \phi \in P_r(T)$$
\end{defi}

For further analysis, several projection operators are also defined in this paper. For cell $T\in T_h$, let $Q_0$ denotes a $L^2$ projection from $[L^2(T)]^d$ to $[P_k(T)]^d$. For each edge  $e\in\epsilon_h$, let $Q_b$ denotes a projection from $[L^2(e)]^d$ to $[P_{k-1}(e)]^d$. We also combine $Q_0$ with $Q_b$, writing $Q_h=\{Q_0,Q_b\}$. Similarly, we also use $\mathbb{Q}_h$ and $\mathbf{Q}_h$ denotes a $L^2$ projection on $P_k(T)$ and $[P_k(T)]^{d\times d}$. And they meet following exchange property, the proof of which can be found in \cite{Wang2016AWG}
\begin{lemma}
Projection $Q_h,\mathbf{Q}_h,\text{and}\mathbb{Q}_h$ meet following exchange property
$$\nabla_w(Q_h \mathbf{v})=\mathbf{Q}_h(\nabla \mathbf{v}),\quad \forall \mathbf{v}\in [H^1(\Omega)]^d,$$
$$\nabla_w\cdot(Q_h \mathbf{v})=\mathbb{Q}_h(\nabla\cdot \mathbf{v}),\quad \forall \mathbf{v}\in H(div,\Omega),$$
\end{lemma}
Next we can define the velocity WG space as follows:
$$V_{h,0}=\{\mathbf{v}=\{\mathbf{v}_0,\mathbf{v}_b\}:\mathbf{v}_0|_{T}\in [P_k(T)]^d,\text{and }\mathbf{v}_b|_e\in [P_{k-1}(T)]^d,\text{ and }\mathbf{v}_b=0\text{ on }\partial \Omega.\}$$
where $P_k(T)$ denotes the space of polynomials of order k or less on T, $P_{k-1}(T)$ denotes the space of polynomials of order k-1 or less on e.
For pressure ,we can give following WG space.
$$Q_h=\{q:q\in L_0^2(\Omega),q|T \in P_{k-1}(T).\}$$
Then, we give WG scheme of (4.2). Define the following WG space:
$$V_h=\{\mathbf{v}\in V_{h,0}|c_w(\mathbf{v},q)=0,\forall q\in Q_h.\}$$
where $c_w(\mathbf{v},q)=(\nabla_w \cdot \mathbf{v},q)$.\\
We will give three bilinear forms defined on $V_h$ as follows. For $\mathbf{w},\mathbf{v}$, define
$$s(\mathbf{w},\mathbf{v})=\gamma(h)\displaystyle\sum_{T\in T_h}h_T^{-1}\left \langle Q_b \mathbf{w}_0-\mathbf{w}_b,Q_b \mathbf{v}_0-\mathbf{v}_b\right \rangle_{\partial T},$$
$$a_w(\mathbf{w},\mathbf{v})=(\nabla_w \mathbf{w},\nabla_w \mathbf{v})+s(\mathbf{w},\mathbf{v}),$$
$$b_w(\mathbf{w},\mathbf{v})=(\mathbf{w}_0,\mathbf{v}_0).$$
where $\gamma(h)$ is taken as follows.
$$\gamma(h)=h^{\epsilon},\quad \epsilon\text{ is a enough small positive number.}$$\\
or
$$\gamma(h)=-\frac{1}{log(h)}.$$
Give that WG scheme (algorithm1)
\begin{algorithm1}
Find $(\lambda_h,\mathbf{u}_h)\in\mathbb{R}\times V_h\text{ such that }b_w(\mathbf{u}_h,\mathbf{u}_h)=1,\text{ and}$ 
\begin{eqnarray}\label{WG-scheme}
 a_w(\mathbf{u}_h,\mathbf{v})=\lambda_h b_w(\mathbf{u}_h,\mathbf{v}),\quad \forall \mathbf{v}\in V_h.
\end{eqnarray}
\end{algorithm1}

\section{error estimate}
In this section, we will verify the assumption presented in section 3, then give the error estimate of WG method by abstract framework in section 3.
Let $V_c=V_0$ and $V=V_c +V_h$. $\forall \mathbf{w},\mathbf{v} \in V$ , we can give an inner product defined as follows:
$$(\mathbf{w},\mathbf{v})_V=(\nabla\mathbf{w}_0,\nabla\mathbf{v}_0)+\displaystyle\sum_{T\in T_h}h_T^{-1}\left \langle Q_b \mathbf{w}_0-\mathbf{w}_b,Q_b \mathbf{v}_0-\mathbf{v}_b\right \rangle_{\partial T}$$
where $\mathbf{w}_0$ denotes  the value of $\mathbf{w}$ inside $T$,$\mathbf{w}_b$ denotes the value of $\mathbf{w}$ on $\epsilon_h$. The corresponding seminorm is defined as follows:

$\parallel \mathbf{w}\parallel_V^2=\displaystyle\sum_{T\in T_h}\parallel \nabla\mathbf{w}_0\parallel_T^2 +\displaystyle\sum_{T\in T_h}\parallel Q_b \mathbf{w}_0-\mathbf{w}_b\parallel_{\partial T}^2$

Obviously, $\parallel \cdot \parallel_V$ coincides $V_c$ with $|\cdot|_1$, and $\parallel\cdot\parallel_V$ defines a norm on $V_h$. Hence $\parallel\cdot\parallel_V$ defines a norm in Hilbert space. Firstly, we verify assumption 1. Obviously, $a(\cdot,\cdot)$ and $a_w(\cdot,\cdot)$ are symmetric bounded linear operator defined on V, and easy to know $a(\cdot,\cdot)$. Hence, we just need to 
prove that $a_w(\cdot,\cdot)$ is positive, referring to 5, we can directly give following result.
\begin{lemma}
$\forall\mathbf{v}\in V_h$, the following inequality holds,
$$a_w(\mathbf{v},\mathbf{v})\gtrsim \gamma(h)\parallel \mathbf{v} \parallel_V^2$$
\end{lemma}

Hence, easy to know assumption 1 holds. Next is assumption 2. we defined following two operator $K$ and $K_h$ in section 3.
To prove $K$ is compact, we define that $$V_1=\{\mathbf{v}\in [H_0^2(\Omega)]^d|c(\mathbf{v},q)=0,\quad\forall q\in Q\}.$$ We can uniquely decompose it into a compound of two operators, i.e. $K=Q\circ K^{\prime}$, where Q is an identity mapping  from $V_1$ to $V_0$, $K^{\prime}$ is an mapping from $V_0$ to $V_1$. On one hand, from classical PDE analysis, when $f\in L^2(\Omega)$, and region $\Omega$ is convex, $\partial \Omega$ is Lipschiz, there exists constant M satisfies $\parallel K^{\prime}f\parallel_{H^2(\Omega)}\leq M \parallel f\parallel_{L^2(\Omega)} $. Hence we know $K^{\prime}$ is bounded operator. We prove $Q$ is compact as follows, we just need to prove for any bounded sequences in $V_1$, both are convergent sequences in $V_0$, we give the following proof.\\
Proof: Since $H^2$ is compact embedded into $H^1$, then $\forall \{\mathbf{v}_n\}_{n=1}^{\infty}\subset V_1$ is bounded, $\exists \mathbf{v}\in [H_0^1(\Omega)]^d, s.t. \mathbf{v}_n\rightarrow v.$ Hence, just need to prove $c(\mathbf{v},q)=0,\forall q\in Q.$
From $\mathbf{v}_n\in [H_0^2(\Omega)]^d$ is bounded, we know $\nabla\cdot \mathbf{v}_n$ is bounded, hence $\displaystyle\lim_{n\to\infty}\nabla\cdot \mathbf{v}_n=\nabla \cdot \mathbf{v}$, and $c(\mathbf{v}_n,q)=0,\forall q\in Q.$ Combing with Lebesgue Dominated convergence theorem, we get $\displaystyle\lim_{n\to\infty}(\mathbf{v}_n,q)=(\mathbf{v},q),\forall q\in Q. $ Hence, operator $Q$ is compact. In summary, operator $K$ is compact.
From the linear, bounded and finite rank properties of operator $K_h$. We can get operator $K_h$ is compact. Hence, assumption 2 holds.\\
We have defined operator $Q_h=\{Q_0,Q_b\}$, where $Q_0$ is $L^2$ projection onto $[P_k(T)]^d$ defined on $T \in T_h$, and $Q_b$ is $L^2$ projection onto $[P_{k-1}(e)]^d$ defined on $e\in\epsilon_h$, then we need to verify $Q_h$ is a projection from $V_0$  to $V_h$, the proof is as follows
:$$\forall\mathbf{u} \in V_0,c_w(Q_h \mathbf{u},q)=(\nabla_w\cdot Q_h \mathbf{u},q) =(\mathbb{Q}_h\cdot\nabla \mathbf{u},q)=(\nabla\cdot \mathbf{u},q)=0.$$Hence, we know $Q_h$ is a projection from $V_0$ to $V_h$.
In addition, $Q_h \mathbf{v}_h=\mathbf{v}_h,\forall \mathbf{v}_h \in V_h$, $Q_0$ is $L^2$ orthogonal projection onto $[P_k(T)]^d$, we have $\forall u\in V$,
$$b(\mathbf{u},\mathbf{v}_h)=(\mathbf{u}_0,\mathbf{v}_0)=(Q_0 \mathbf{u}_0,\mathbf{v}_0)=b(Q_h \mathbf{u},\mathbf{v}_h).$$
Hence assumption 3 holds.
Next, we will verify assumption 4 and 5 respectively, and give an estimate of eigenvector based on those two assumptions.\\Let X $[L^2(\Omega)]^d$. Denote
$$\delta_{h,\mu}=\displaystyle\sup_{\mathop{}_{\parallel \mathbf{v}\parallel_V=1}^{\mathbf{v}\in R(E_{\mu}(K))}}\parallel \mathbf{v}-Q_h \mathbf{v}\parallel_V,$$
$$\delta_{h,\mu}^{\prime}=\displaystyle\sup_{\mathop{}_{\parallel \mathbf{v}\parallel_V=1}^{\mathbf{v}\in R(E_{\mu}(K))}}\parallel \mathbf{v}-Q_h \mathbf{v}\parallel_X,$$
For defined errors $\delta_{h,\mu}$ and $\delta_{h,\mu}^{\prime}$, the following estimate holds, the method of proof can be found in the \cite{zhai2019weak}.
\begin{lemma}
Denote $R(E_{\mu}(K))\in [H^k(\Omega)]^d$, we have that
$$\delta_{h,\mu}\lesssim h^k$$
$$\delta_{h,\mu}^{\prime}\lesssim h^{k+1}.$$
\end{lemma}

Denote
$$e_{h,\mu}=\parallel (K-K_h Q_h)|_{R(E_{\mu}(K))})\parallel_V,$$
$$e_{h,\mu}^{\prime}=\parallel (\Pi_0 K-\Pi_0 K_h Q_h)|_{R(E_{\mu}(K))})\parallel_X.$$
We will estimate $e_{h,\mu}$ and $e_{h,\mu}^{\prime}$ by similar method to \cite{Xie2015TheWG} . Firstly, we give the equivalent variational form of Stokes equation and Stokes WG scheme as follows:
\begin{equation}
    \begin{cases}
        -\Delta \mathbf{u}+\nabla p=\mathbf{f},& \text{in }\Omega ,\\
        \nabla \cdot \mathbf{u}=0,& \text{in }\Omega,\\ 
        \mathbf{u}=\mathbf{0},& \text{on }\partial\Omega , \\
     \end{cases}
\end{equation}
where $\mathbf{f}\in [L^2(\Omega)]^d$. 
The corresponding variational form: Find $(\mathbf{u},p) \in  [H_0^1(\Omega)]^d \times Q$ such that\\
\begin{equation}
    \begin{cases}
        a(\mathbf{u},\mathbf{v})-c(\mathbf{v},p)=(\mathbf{f},\mathbf{v}), \forall \mathbf{v} \in [H_0^1(\Omega)]^d,\\
        c(\mathbf{u},q)=0,      \forall q \in Q,
     \end{cases}
\end{equation}\\
where $a(\cdot,\cdot)$ and $c(\cdot,\cdot)$ agree with former definition.
\ The corresponding WG scheme: Find $(\mathbf{u}_h,p_h) \in  V_{h,0} \times Q_h$ such that \\
\begin{equation}
    \begin{cases}
        a_w(\mathbf{u}_h,\mathbf{v}_h)-c_w(\mathbf{v}_h,p_h)=(\mathbf{f},\mathbf{v}_0), \forall \mathbf{v}_h \in V_{0,h},\\
        c_w(\mathbf{u}_h,q_h)=0,      \forall q_h \in Q_h,
     \end{cases}
\end{equation}\\
where $a_w(\cdot,\cdot)$ and $c_w(\cdot,\cdot)$ agree with former definition.
Further, from (3.7) and (3.8) we can get following variational problem (3.9) and (3.10).\\
The corresponding variational form: Find $\mathbf{u} \in  V_0 $ such that \\
\begin{equation}
        a(\mathbf{u},\mathbf{v})=(\mathbf{f},\mathbf{v}), \forall \mathbf{v} \in V_0,\\
\end{equation}\\
The corresponding WG scheme:Find $\mathbf{u}_h \in  V_h $ such that \\
\begin{equation}
        a_w(\mathbf{u}_h,\mathbf{v})=(\mathbf{f},\mathbf{v}_0), \forall \mathbf{v} \in V_h,\\
\end{equation}\\
From the proof of the compactness of $K$ we can know $V_0$ is a closed subspace of  $[H_0^1(\Omega)]^d$, similarly, we can prove $V_h$ is closed sunset of $V_{h,0}$, and $V_0$ and $V_h$ both are Hilbert space, and from the proof of assumption 1 we know $a(\cdot,\cdot)$ and $a_w(\cdot,\cdot)$ are bounded, positive bilinear forms defined on $V_0$ and $V_h$, respectively. Hence from Lax-Milgram theorem, problem (5.4) and (5.5) both have unique solution, hence easy to know (5.2) and (5.4) is equivalent, and (5.3) and (5.5) is equivalent.

Notice two operator $K:V_c\rightarrow V_c$ and $K_h:V_h\rightarrow V_h$ defined former, they meet the following conditions :
$$a(K\mathbf{f},\mathbf{v})=b(\mathbf{f},\mathbf{v}),\quad \forall \mathbf{v}\in V_c,$$
$$a_w(K_h \mathbf{f}_h,\mathbf{v}_h)=b(\mathbf{f}_h,\mathbf{v}_h),\quad \forall \mathbf{v}_h\in V_h.$$
Next, we will show $\mathbf{u}:=K\mathbf{f},\mathbf{u}_h:=K_h Q_h \mathbf{f},\forall \mathbf{f}\in [L^2(\Omega)]^d$ are the exact solution and numerical solution of (5.1) respectively.
Firstly, for $K$, $a(K\mathbf{f},\mathbf{v})=b(\mathbf{f},\mathbf{v}),\forall \mathbf{v}\in V_0,$ let $u:=Kf$, then,
$$a(\mathbf{u},\mathbf{v})=(\mathbf{f},\mathbf{v}),\forall \mathbf{v}\in V_0,$$
From the equivalence of variational problem, $\mathbf{u}:=K\mathbf{f}$ is the exact solution of (5.1).
Similarly, for $K_h$, we have $a_w(K_h \mathbf{f}_h,\mathbf{v})=b(\mathbf{f}_h,\mathbf{v}),\forall \mathbf{v}\in V_h,$
Denote $\mathbf{f}_h=Q_h \mathbf{f},\mathbf{u}_h=K_h Q_h \mathbf{f}$, then $$a_w(\mathbf{u}_h,\mathbf{v})=b(Q_h \mathbf{f},\mathbf{v})=b(\mathbf{f},\mathbf{v})=(\mathbf{f},\mathbf{v}_0),\forall \mathbf{v}\in V_h,$$
From the equivalence of variational problem, $\mathbf{u}_h:=K_h Q_h \mathbf{f}$ is the numerical solution of (5.1).

Hence, we can give an estimate of $e_{h,\mu}$ and $e_{h,\mu}^{\prime}$ by above analysis and following error estimate of Stokes equation. As for the error estimate of Stokes equation, we just give a conclusion and the concrete process can be found in appendix A.
\begin{theorem}
Denote $\mathbf{u}\in [H^{k+1}(\Omega)]^d$ exact solution of (5.2), $\mathbf{u}_h$ is numerical solution of (5.3), and dual problem of (5.2) has $[H^2(\Omega)]^d\times H^1(\Omega)$-regularity. Then we have following estimate,
$$\parallel \mathbf{u}-\mathbf{u}_h\parallel_V\leq C\gamma(h)^{-1}h^k(\parallel \mathbf{u}\parallel_{k+1}+\parallel p\parallel_k),$$
$$\parallel \mathbf{u}-\mathbf{u}_h\parallel_X\leq C\gamma(h)^{-1}h^{k+1}(\parallel \mathbf{u}\parallel_{k+1}+\parallel p\parallel_k).$$
\end{theorem}
\begin{lemma}
 If $R(E_{\mu}(K))\subset [H^{k+1}(\Omega)]^d$, then
 $$e_{h,\mu}\lesssim\gamma(h)^{-1}h^k,$$
 $$e_{h,\mu}^{\prime}\lesssim \gamma(h)^{-1}h^{k+1}.$$
 \end{lemma}

 Proof: from theorem 5.3, we know \\$e_{h,\mu}=\parallel(K-K_h Q_h)|_{R(E_{\mu}(K))}\parallel_V$
 $=\displaystyle\sup_{\parallel \mathbf{f} \parallel_V=1}\parallel K\mathbf{f}-K_h Q_h \mathbf{f}\parallel_V$\\
 $=\displaystyle\sup_{\parallel \mathbf{f}\parallel_V=1}\parallel \mathbf{u}-\mathbf{u}_h\parallel_V$
 $\lesssim \gamma(h)^{-1}h^k.$\\
 $e_{h,\mu}^{\prime}=\parallel(\Pi_0 K-\Pi_0 K_h Q_h)|_{R(E_{\mu}(K))}\parallel_X$
 $=\displaystyle\sup_{\parallel \mathbf{f} \parallel_X=1}\parallel \Pi_0 K\mathbf{f}-\Pi_0 K_h Q_h \mathbf{f}\parallel_X$\\
 $=\displaystyle\sup_{\parallel \mathbf{f}\parallel_X=1}\parallel \Pi_0 (\mathbf{u}-\mathbf{u}_h)\parallel_X$
  $\leq\displaystyle\sup_{\parallel \mathbf{f}\parallel_X=1}\parallel  \mathbf{u}-\mathbf{u}_h\parallel_X$
 $\lesssim \gamma(h)^{-1}h^{k+1}.$
 
 Hence, from theorem 3.1 and 3.2, we have following estimate of eigenfunction.
 \begin{theorem}
 Denote $\lambda$ eigenvalue of problem (5.1) with m multiplicity m, $\{\lambda_{h,j}\}_{j=1}^{\infty}$ is the corresponding numerical eigenvalues. And $R(E_{\mu}(K))\subset [H^{k+1}(\Omega)]^d$ is a m-dimensional eigenspace related to $\lambda$, $\{\mathbf{u}_{j,h}\}_{j=1}^{\infty}$ is a set of bases of eigenspace $R(E_{\mu}(K))$ related to $\{\lambda_{h,j}\}_{j=1}^{\infty}$. Then, when h is small  enough, there exists a set of eigenfunctions $u_j\in R(E_{\mu}(K)),j=1,\cdots,m$ such that
 $$\parallel \mathbf{u}_j-\mathbf{u}_{j,h}\parallel_ V\lesssim \gamma(h)^{-1}h^k,$$
 $$\parallel \mathbf{u}_j-\mathbf{u}_{j,h}\parallel_ X\lesssim \gamma(h)^{-1}h^{k+1}.$$
  
 where 
 $$\gamma(h)=h^{\epsilon} \text{ for small enough positive number }\epsilon,$$
 or
 $$\gamma(h)=-\frac{1}{log(h)}.$$
 Hence, we have following result
 $$\parallel \mathbf{u}_j-\mathbf{u}_{j,h}\parallel_ V\lesssim h^{k-\epsilon},$$
 $$\parallel \mathbf{u}_j-\mathbf{u}_{j,h}\parallel_ X\lesssim h^{k+1-\epsilon}.$$
 or
  $$\parallel \mathbf{u}_j-\mathbf{u}_{j,h}\parallel_ V\lesssim -log(h)h^k,$$
 $$\parallel \mathbf{u}_j-\mathbf{u}_{j,h}\parallel_ X\lesssim -log(h)h^{k+1}.$$
 \end{theorem}
 Next we will verify assumption 6 and give an estimate of eigenvalue based on this. Firstly, we need estimate $\epsilon_{h,\mathbf{u}}=a(\mathbf{u},\mathbf{u})-a_w(Q_h \mathbf{u},Q_h  \mathbf{u})$.
 \begin{lemma}
 $\forall\mathbf{u}\in [H^{k+1}(\Omega)]^d$, then,
 $$|\epsilon_{h,\mathbf{u}}|\lesssim h^{2k}.$$
 \end{lemma}
 
 From theorem 3.9, we can get the error estimate of eigenvalue.
 \begin{theorem}
  Denote $\lambda$ eigenvalue of (5.1) with m multiplicity, $\{\lambda_{h,j}\}_{j=1}^{\infty}$ is the corresponding numerical eigenvalue, $R(E_{\mu}(K))\subset [H^{k+1}(\Omega)]^d$ is m dimension eigenspace related to $\lambda$, $\{\mathbf{u}_{j,h}\}_{j=1}^{\infty}$ is a set of bases of eigenspace $R(E_{\mu}(K))$ related to  $\{\lambda_{h,j}\}_{j=1}^{\infty}$. Then $\forall\ 1\leq j\leq m$, when h is small enough, the following estimate holds,
$$|\lambda-\lambda_{h,j}|\lesssim h^{2k}.$$
 \end{theorem}
 Proof: from theorem 3.3, we have 
 $$|\lambda-\lambda_{h,j}|\leq C(\epsilon_{h,\mathbf{u}_j}++e_{h,\lambda^{-1}}^2+e_{h,\lambda^{-1}}^{\prime  2}+\delta_h^2\gamma(h)^{-2}+\delta_h^{\prime2}\gamma(h)^{-2})$$ and from lemma 5.2,5.4,5.6 we have
$\epsilon_{h,\mathbf{u}_j}\lesssim h^{2k},\quad \delta_{h,\mu}\lesssim h^k,\quad\delta_{h,\mu}^{\prime}\lesssim h^{k+1},$\\
$$e_{h,\mu}\lesssim \gamma(h)^{-1}h^k,\quad e_{h,\mu}^{\prime}\lesssim \gamma(h)^{-1}h^{k+1},$$\\
This shows 
$$|\lambda-\lambda_{h,j}|\lesssim\gamma(h)^{-2}h^{2k}.$$
 In fact, all estimates we do so far are no dependent on the choice of $\gamma(h)$, we can just take $\gamma(h)=1$ to get best estimate,
 $$|\lambda-\lambda_{h,j}|\lesssim h^{2k}$$

\section{asymptotic lower bounds}
In this section we prove WG scheme can give asymptotic lower bounds of the eigenvalues.
 From the abstract framework in section three, we know we can establish the lower bounds
by verifying assumption seven. To verify assumption 7, we introduce a lower bounded estimate from \cite{Lin2014LowerBO}, this estimate plays a key role in the following proof.

 Denote $\mathbf{u}$ the exact solution of eigenvalue problem (5.1), then the following estimate holds:
 $$\displaystyle\sum_{T\in T_h}\parallel \nabla \mathbf{u}-Q_h \nabla \mathbf{u}\parallel_T^2\gtrsim Ch^{2k}.$$
 \begin{lemma}
 Suppose $(\lambda,\mathbf{u})$ is exact solution of (5.1), and $(\lambda_h,\mathbf{u}_h)$ is the corresponding numerical solution. Denote $\gamma(h)\ll 1$, then for small enough h, the following equation holds,
 $$\epsilon_{h,\mathbf{u}}\geq \lambda_h \parallel \mathbf{u}-\mathbf{u}_h\parallel_X^2.$$
 
 \end{lemma}
 Proof:
 $\epsilon_{h,\mathbf{u}}=a(\mathbf{u},\mathbf{u})-a_w(Q_h \mathbf{u},Q_h \mathbf{u})$
 \\$=\parallel \nabla \mathbf{u}\parallel^2-\displaystyle\sum_{T\in T_h}\parallel\nabla_w Q_h \mathbf{u}\parallel_T^2-s(Q_h \mathbf{u},Q_h \mathbf{u})$
 \\$=\parallel \nabla \mathbf{u}\parallel^2-\displaystyle\sum_{T\in T_h}\parallel\nabla_w Q_h \mathbf{u}\parallel_T^2-\displaystyle\sum_{T\in T_h}\gamma(h)h_T^{-1}\parallel Q_b(Q_0 \mathbf{u}-Q_b \mathbf{u})\parallel_{\partial T}^2$\\
 and $\displaystyle\sum_{T\in T_h}\parallel \nabla \mathbf{u}-\mathbf{Q}_h \nabla \mathbf{u}\parallel_T^2\gtrsim h^{2k},$ and $\displaystyle\sum_{T\in T_h}\gamma(h)h_T^{-1}\parallel Q_b(Q_0 \mathbf{u}-Q_b \mathbf{u})\parallel_{\partial T}^2\lesssim \gamma(h)h^{2k}$\\
 since $\gamma(h)\ll 1$, then when h is small enough, $$\epsilon_{h,\mathbf{u}}\gtrsim h^{2k}.$$
And from theorem 3.4 we know $\lambda_h \parallel \mathbf{u}-\mathbf{u}_h\parallel_X^2\lesssim h^{2k+2}.$
 Hence, for h is small enough, the following equation holds,
 $$\epsilon_{h,\mathbf{u}}\geq \lambda_h\parallel \mathbf{u}-\mathbf{u}_h\parallel_X^2,$$
 \begin{theorem}
Suppose $(\lambda,\mathbf{u})$ is exact solution of (5.1) and $(\lambda_h,\mathbf{u}_h)$ is the corresponding numerical solution. And when $h\rightarrow0, \gamma(h)\rightarrow0$, for small enough $h$, the following equation holds,
$$\lambda\geq\lambda_h.$$
 \end{theorem}
 Proof: from lemma 6.1 we know
 $$\displaystyle\sum_{T \in T_h}\parallel \nabla \mathbf{u}-\mathbf{Q}_h\nabla \mathbf{u} \parallel_T^2\gtrsim h^{2k}$$
 $$\displaystyle\sum_{T \in T_h}h_T^{-1}\parallel Q_b(Q_0 \mathbf{u}-Q_b \mathbf{u})\parallel_{\partial T}^2\lesssim \gamma(h) h^{2k}.$$
 Hence, there exists constant $C_0$ independent of $h$, such that
 $$\displaystyle\sum_{T \in T_h}h_T^{-1}\parallel Q_b(Q_0 \mathbf{u}-Q_b \mathbf{u})\parallel_{\partial T}^2\leq C_0 \displaystyle\sum_{T\in T_h}\parallel \nabla \mathbf{u}- \mathbf{Q}_h \nabla \mathbf{u}\parallel_T^2$$

\section{guarantee lower bounded}

In section 6, we prove the WG scheme (5.1) provides asymptotic lower bounds of
the eigenvalues, but we can’t prove the lower bounds property of a constant in this way, so we will take a more sophisticated method to prove the guaranteed error control of the stokes eigenvalues referring to \cite{carstensen2020skeletal}.
Firstly, we introduce several inequality on $H_0^1(\Omega)$ and a more complex stabilizer provided by \cite{carstensen2020skeletal}.
 $$\parallel (1-\mathbb{Q}_0)f\parallel \leq \delta\parallel(1-Q)\nabla f\parallel^2\quad \forall f\in H_0^1(\Omega).$$
 $$s(\mathbb{Q}_h f,\mathbb{Q}_h f)\leq \alpha \Lambda\parallel (1-Q)\nabla f\parallel \forall f\in H_0^1(\Omega).$$
 $$s(\mathbf{u}_h, \mathbf{v}_h) := \frac{\alpha}{n+1}\sum_{T\in \tau} \sum_{F\in \mathcal{F}(T)} h_T^{-2}
| F|^{-1} |T| \langle Q_b \mathbf{u}_0-\mathbf{u}_b, Q_b \mathbf{v}_0-\mathbf{v}_b\rangle_{[L^2(F)]^d},$$

 And \cite{carstensen2020skeletal} provides the theorem to give the conditions for those inequalities to hold as follows:

\begin{theorem}
For space $V_{h,0}$ and $Q_h$ and suppose $C_{apx}$ satisfies the following inequality. Then (A)-(B) hold with $\delta=C_{apx}^2 h_{max}^2$ and $\Lambda=(C_{apx}+\frac{2}{n+1})\frac{C_{apx}}{n+1}$,
$$\parallel f-\mathbf{\Pi}_k f\parallel_{L^2(\Omega)}\leq C_{apx} h_T \parallel (1-\mathbf{\Pi}_{k-1}\nabla f)\parallel_{L^2(\Omega)}\quad \forall f\in H^1(\Omega).$$
where $\mathbf{\Pi}_k$ is the $L^2$ projection onto $P_k(\Omega).$
\end{theorem}

Since $Q_0$ is a projection from $V_0$ to $V_h$, we can consider $Q_0$ as the restriction of $Q_0$ defined on $[H_0^1(\Omega)]^d$ on $V_0$. Hence, by similar proof method, we can give those inequalities and its conditions to hold on $V_0$ as follows.
\begin{theorem}
For space $V_{h,1}:=\{v=\{v_0,v_b\}: v_0|_{T}\in P_k(T),\text{ and }v_b|_e\in P_{k-1}(T),\text{ and }\mathbf{v}_b=0\text{ on }\partial T\}$ and $Q_h$ and suppose $C_{apx}$ satisfies the following inequality. Then (A)-(B) hold with $\delta=C_{apx}^2 h_{max}^2$ and $\Lambda=(C_{apx}+\frac{2}{n+1})\frac{C_{apx}}{n+_1}$,
$$\parallel \mathbf{f}-\Pi_k \mathbf{f}\parallel_{L^2(\Omega)}\leq C_{apx} h_T \parallel (1-\Pi_{k-1}\nabla \mathbf{f})\parallel_{L^2(\Omega)}\quad for\ all f\in H^1(\Omega).$$
where $\Pi_k$ is the $L^2$ projection onto $[P_k(\Omega)]^d.$
\end{theorem}

Based on above inequality, we can use exactly the same method to prove the GLB property of Stokes eigenvalue problem. Hence, we can get the following theorem,
\begin{theorem}
If $\delta$  and $\Lambda$  in (A)-(B) and the stability parameter $\alpha\geq 0$,
satisfy either (i) $\delta \lambda  + \alpha \Lambda  \leq$  1 or (ii) $\delta \lambda h + \alpha \Lambda  \leq  1$, then $\lambda_h \leq  \lambda$.
\end{theorem}

Next, we will consider the lowest-order case (1,0,0). Firstly, we defined a new operator $I_C=\{I_{NC},\Pi_0\}$, where $I_{NC}\mathbf{f}\in [P_1(\Omega)]^d$, and $I_{NC} (\mathbf{\phi})(mid(F))_i:=\frac{\int_F \mathbf{\phi_i} ds}{|F|}, i=1,\cdots,d.$ We need to show $I_C$ has similar exchange properties to $Q_h$, hence $I_{C}$ defined on $V_0$ in can be seen the restriction of $I_{C}$ defined on $[H_0^1(\Omega)]^d$ on $V_0$. The proof is as follows:
$$\forall\mathbf{u} \in V_0,c_w(I_c \mathbf{u},q)=c_w(Q_h \mathbf{u},q)=(\mathbb{Q}_h\cdot\nabla \mathbf{u},q)=0.$$
Based on just-defined operator $I_{NC}$, we can use exactly the same method to prove the GLB property of Stokes eigenvalue problem of the lowest-order situation. Hence, we can get the following theorem.
\begin{theorem}
If the lowest-order situation with maximal mesh-size $h_{max}$ and stability parameter $\alpha \geq 0$ satisfies $max\bigl\{\alpha, min\{ \lambda , \lambda_h\}  h_{max}^2\bigr\} 
\leq  \kappa_{CR}^{-2}.$ for the m-th
discrete eigenvalue $\lambda_h$, then $\lambda_h$ is a GLB For the m-th exact eigenvalue $\lambda  \geq  \lambda_h.$

\end{theorem}

\section{numerical results}
In this section, we will show some numerical results of WG scheme analyzed in former sections.
\subsection{Stokes eigenvalue problem in unit square region}
In this example, we consider problem (5.1) on the square region $\Omega = (0,1)^2$. Since the exact eigenvalues of the Stokes eigenvalue problem (5.1) are unknown, we will use the values of the first six eigenvalues of problem (5.1) provided in the article \cite{Xie2018ExplicitLB} as the exact eigenvalues.

In this example, we take uniform triangulation and denote the mesh size by $h$. Take the order of polynomial in WG space is $k=1$, and take parameter $\gamma(h) = 1, h^{0.1}, -\frac{1}{log(h)}$, respectively. The numerical results of the first six eigenvalues are shown in Figure 8.1(a)-(c). From these images, we can find that the errors of the eigenvalues are second-order convergent in agreement with the conclusion in Theorem 5.7. Also, for the case with the above mesh with polynomial order $k=1$ and parameter $\gamma(h)=h^{0.1}$, we present the numerical results for the first six eigenvalues in Table 8.1. Table 8.1 shows that all numerical solutions of the eigenvalues are lower bounds of the exact eigenvalues, which is consistent with Theorem 6.2. Also, the convergence order coincides with Theorem 5.7.

Besides, since the eigenvector of problem (3.1) is unknown and it’s no easy to find a standard value to be a reference for the exact eigenvector. Hence, we can’t give the error estimate of eigenfunction there, but the numerical solution $\lambda=52.344691168$ of the first eigenvalue for case that mesh size $h=\frac{1}{64}$, polynomial order $k=1$ and parameter $\gamma(h)=C$ is given in Figures 8.2(a)(b) and 8.3(a)(b) $ (\ mathbf{u}=\{u_1,u_2\},p)$ in the image, where Figures 8.2, 8.3(a)(b) represent $u_1,u_2,\mathbf{u},p$, respectively.

\begin{table}[htb]
\begin{center}
\caption{error of eigenvalue in rectangle region }  
\label{table:1} 
\resizebox{\linewidth}{!}{
\begin{tabular}{|c|c|c|c|c|c|c|}   
\hline   $h$  & $1/4$& $1/8$ & $1/16$ & $1/32$ & $1/64$ & $1/128$  \\   
\hline    $\lambda_1 -\lambda_{1,h}$ &  2.5122e+1 & 1.0208e+1 & 3.1731e+0 & 8.8337e-1 & 2.3766e-1 & 6.3365e-2  \\ 
\hline   order&  0.7215 & 1.2992& 1.6857 &1.8448  & 1.8940 & 1.9071  \\  
\hline    $\lambda_2 -\lambda_{2,h}$ &  6.0736e+1 & 2.9330e+1 & 9.9479e+0 & 2.8491e+0 & 7.7125e-1 & 2.0553e-1  \\ 
\hline   order  & 0.4189& 1.0501 & 1.5599& 1.8038 & 1.8852 & 1.9078  \\  
\hline    $\lambda_3 -\lambda_{3,h}$ &  5.9361e+1& 2.8763e+1 & 9.7586e+0& 2.7955e+0 & 7.5732e-1 & 2.0201e-1  \\ \hline   order   &  0.4392& 1.0452 & 1.5595&  1.8035 & 1.8841 & 1.9064  \\  
\hline    $\lambda_4 -\lambda_{4,h}$  & 9.1862e+1& 4.8975e+1 & 1.7846e+1& 5.2640e+0 & 1.4395e+0 & 3.8517e-1\\ 
\hline   order  & 0.3417& 0.9074 & 1.4564& 1.7614 & 1.8705 & 1.9020  \\  
\hline    $\lambda_5 -\lambda_{5,h}$  & 1.1734e+2& 6.8730e+1 & 2.6164e+1& 7.8362e+0 & 2.1498e+0 & 5.7496e-1  \\ 
\hline   order  &  0.2742& 0.7716 & 1.3933& 1.7394 & 1.8659 & 1.9026  \\  
\hline    $\lambda_6 -\lambda_{6,h}$ & 1.2931e+2& 7.3735e+1 & 2.8333e+1& 8.5546e+0 & 2.3597e+0 & 6.3381e-1  \\ 
\hline   order  & 0.2202& 0.8104 & 1.3798& 1.7277 & 1.8580 & 1.8964  \\  
 
\hline
\end{tabular}   
}
\end{center}   
\end{table}

\begin{figure}[H]
\centering
\subfigure[$\gamma(h)=C$]{\includegraphics[scale=0.4]{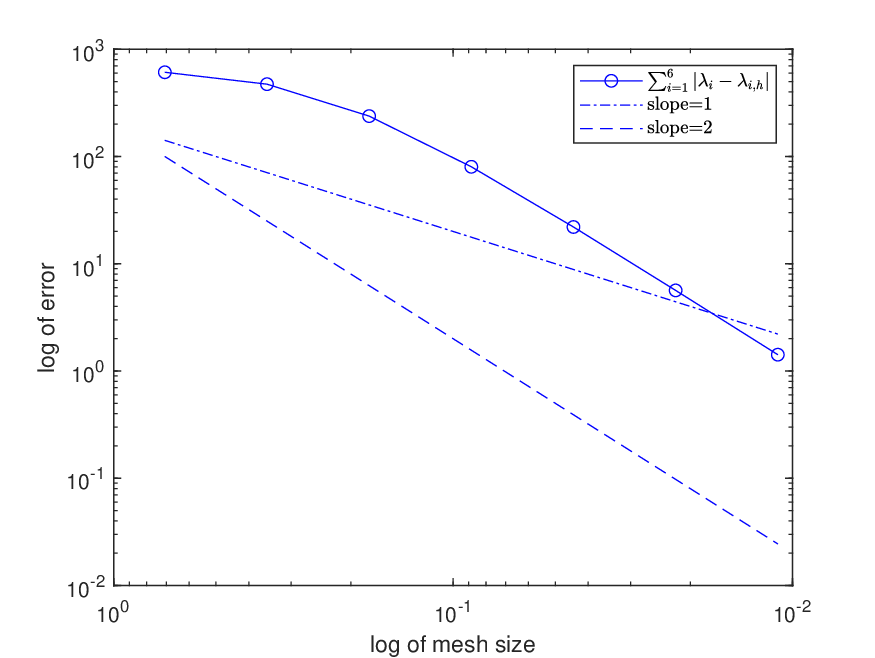}}
\subfigure[$\gamma(h)=h^{0.1}$]{\includegraphics[scale=0.4]{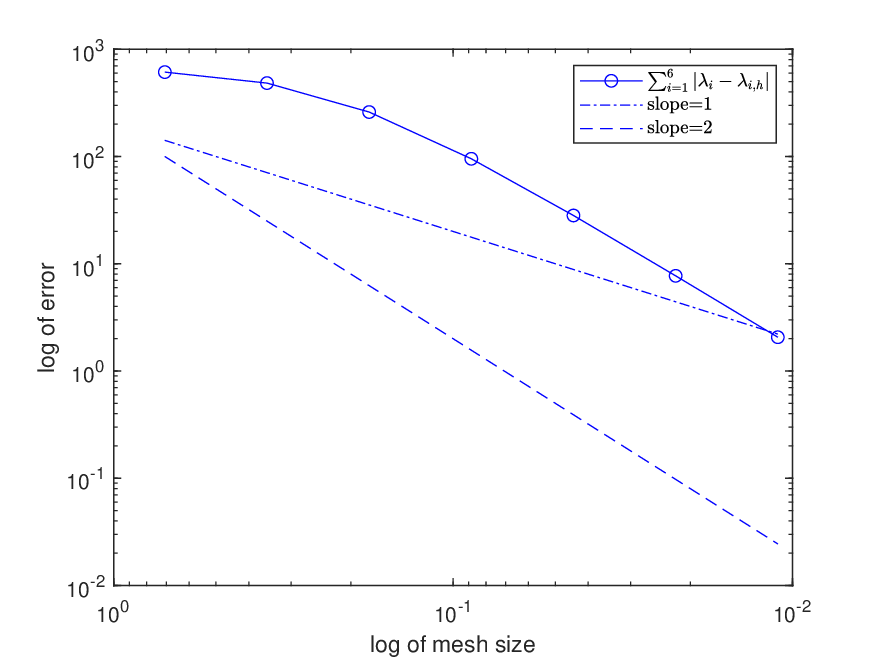}}
\\
\subfigure[$\gamma(h)=-\frac{1}{log(h)}$]{\includegraphics[scale=0.4]{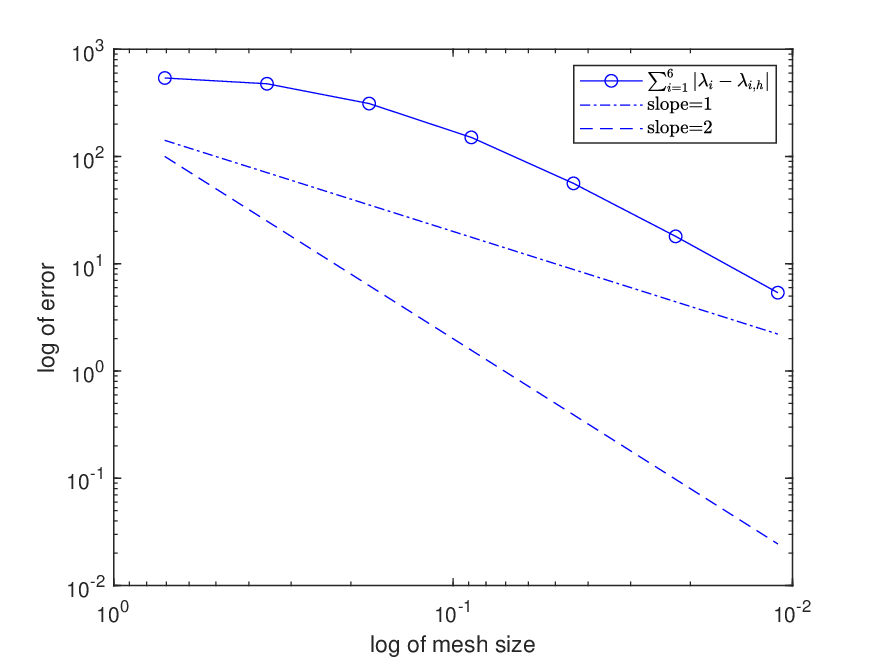}}
\caption{error estimate of eigenvalue in rectangle region}
\end{figure}

\begin{figure}[H]
\centering
\subfigure[]{\includegraphics[scale=0.40]{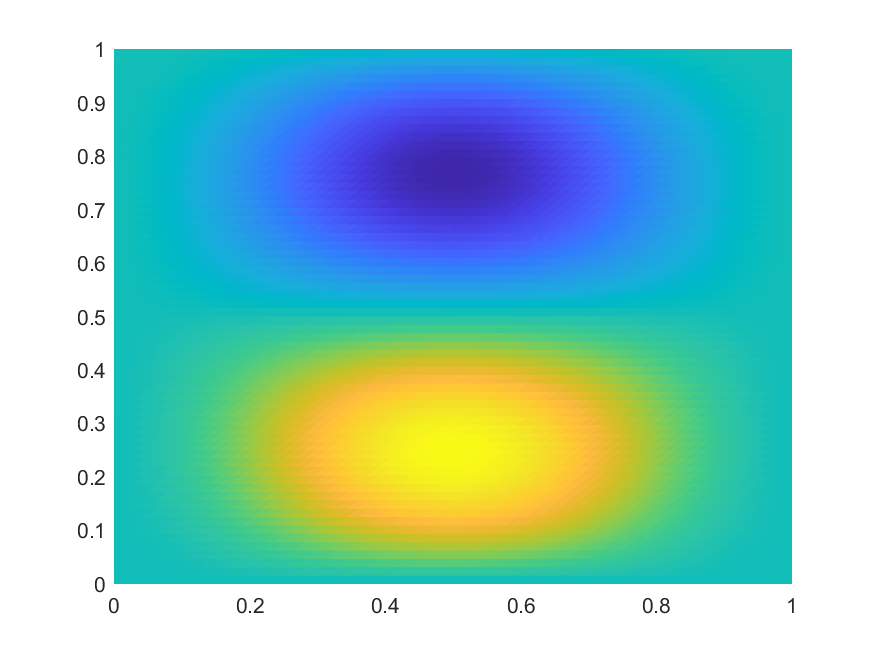}}
\subfigure[]{\includegraphics[scale=0.40]{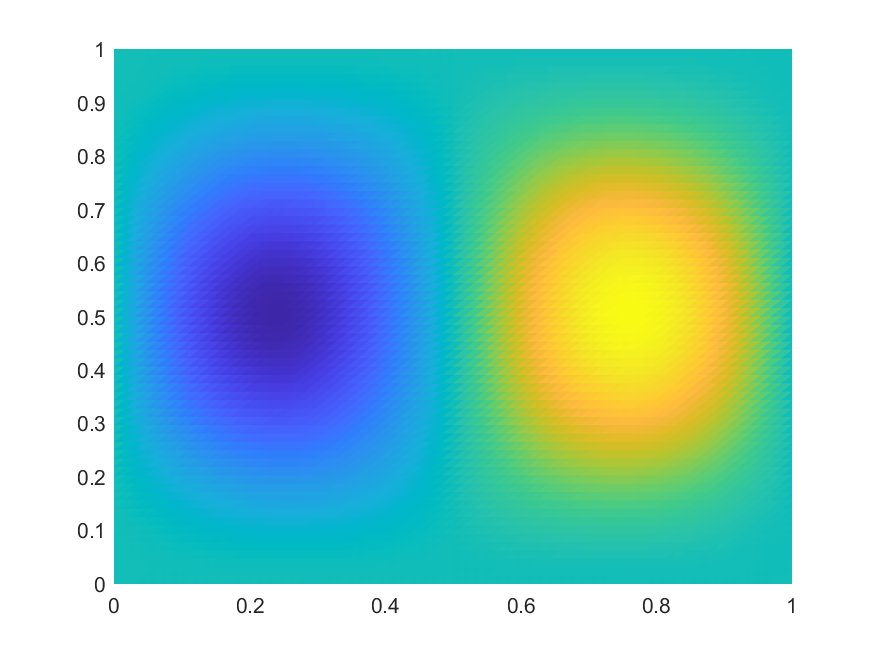}}
\caption{$(a)\text{eigenfunction } \mathbf{u}_1,(b)\text{ eigenfunction } \mathbf{u}_2$}
\end{figure}

\begin{figure}[H]
\centering
\subfigure[]{\includegraphics[scale=0.40]{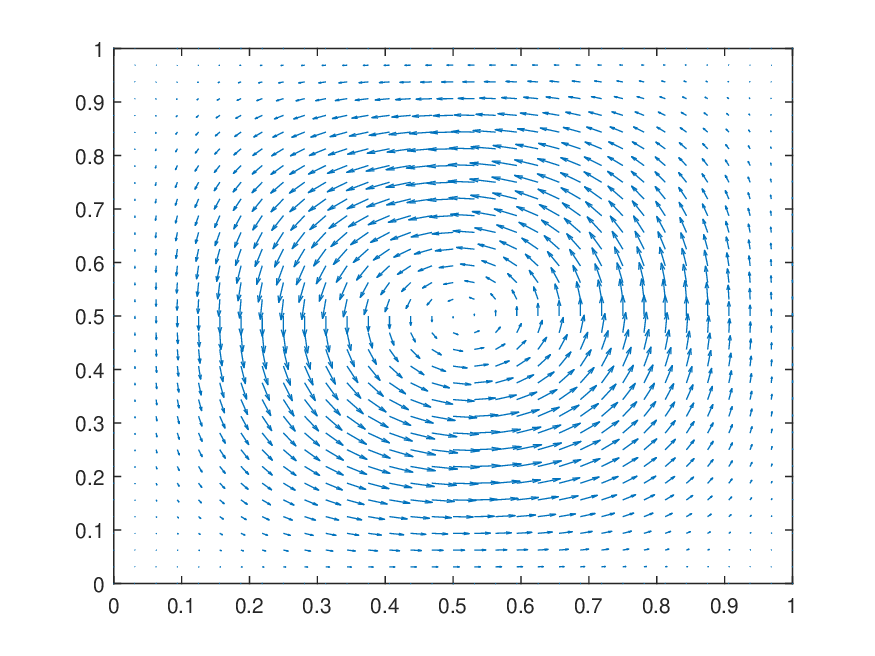}}
\subfigure[]{\includegraphics[scale=0.40]{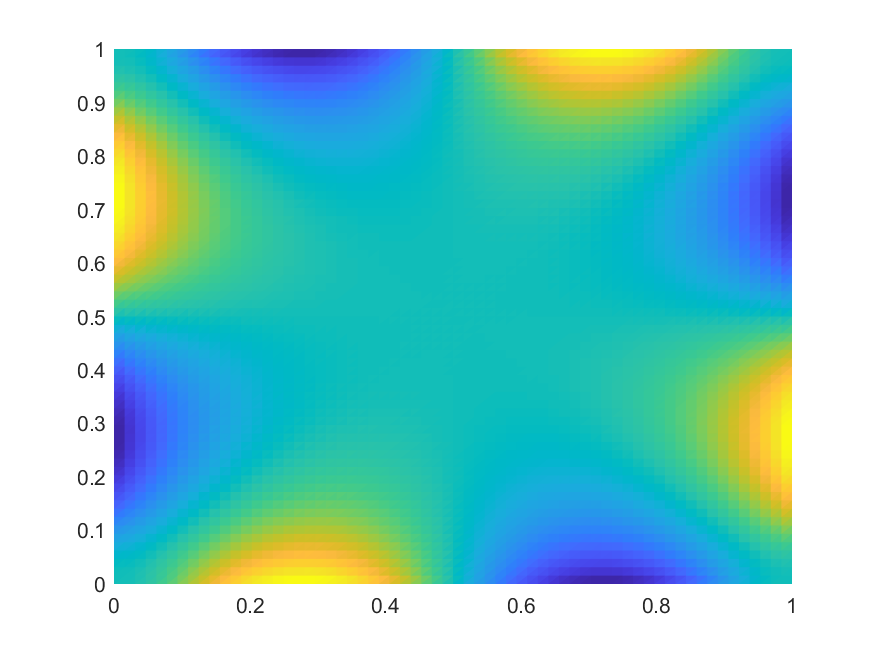}}
\caption{$(a)\text{ eigenfunction }\mathbf{u},(b)\text{ eigenfunction }p$}
\end{figure}

\subsection{ Stokes eigenvalue problem in L-shaped region }
In this example, we consider problem (5.1) on the L-shaped region $\Omega=(-1,1)^2\backslash[0,1]^2$. Since the exact eigenvalues of the Stokes eigenvalue problem (5.1) are unknown, we will use the values of the first five eigenvalues of problem (5.1) provided in the article \cite{Xie2018ExplicitLB} as the exact eigenvalues.

In this example, we take uniform triangulation and denote the mesh size by $h$. Take the order of polynomial in WG space is $k=1$, and take parameter $\gamma(h) = 1, h^{0.1}, -\frac{1}{log(h)}$, respectively. The numerical results of the first five eigenvalues are shown in Figure 8.4(a)-(c). From these images, we can find that the errors of the eigenvalues are second-order convergent in agreement with the conclusion in Theorem 5.7. Also, for the case with the above mesh with polynomial order $k=1$ and parameter $\gamma(h)=h^{0.1}$, we present the numerical results for the first five eigenvalues in Table 8.2. Table 8.2 shows that all numerical solutions of the eigenvalues are lower bounds of the exact eigenvalues, which is consistent with Theorem 6.2. Also, the convergence order coincides with Theorem 5.7.(Note: In Table 8.2 we can notice that for the smallest eigenvalue on the L-shaped region, the convergence order is far from 1.9 or 2, but is around 1.4 and 1.5. This is because the convergence order of the eigenvalue decreases by about 0.544 under uniform mesh grid, due to the nature of the eigenfunction, which is independent of the numerical solution method we use for the eigenvalue problem. More details can be found in \cite{gedicke2018arnold}.)

Besides, since the eigenvector of problem (3.1) is unknown and it’s no easy to find a standard value to be a reference for the exact eigenvector. Hence, we can’t give the error estimate of eigenfunction there, but the numerical solution $\lambda=32.13269465$ of the first eigenvalue for case that mesh size $h=\frac{1}{64}$, polynomial order $k=1$ and parameter $\gamma(h)=C$ is given in Figures 8.5(a)(b) and 8.6(a)(b)(c) $(\mathbf{u}=\{u_1,u_2\},p)$ in the image, where Figures 8.5(a), 8.5(b), 8.6(a)(b), 8.6(c) represent $u_1,u_2,\mathbf{u},p$, respectively. (Note: In Figure 3.6(a), the distribution of the edge region cannot be shown in the figure because the difference between the value of its central region and the value of the edge region (lower right) (upper left) is too large. Therefore, we draw the distribution of the lower right part separately in 8.6(b).)

\begin{table}[H]   
\begin{center}   
\caption{ error of eigenvalue in L-shaped region }  
\label{table:2} 
\resizebox{\linewidth}{!}{
\begin{tabular}{|c|c|c|c|c|c|c|}   
\hline   $h$  & $1/4$& $1/8$ & $1/16$ & $1/32$ & $1/64$ & $1/128$  \\   
\hline    $\lambda_1 -\lambda_{1,h}$ &  2.1729e+1 & 1.2854e+1 &5.1809e+0 & 1.7829e+0 & 6.0681e-1 & 2.2006e-1  \\ 
\hline   order& &  0.7573 & 1.3109& 1.5389 &1.5549  & 1.4633   \\  
\hline    $\lambda_2 -\lambda_{2,h}$ &  2.6304e+1 & 1.6177e+1 & 6.2785e+0 & 1.9308e+0 & 5.4509e-1 & 1.5018e-1  \\ 
\hline   order  &  & 0.7013 & 1.3655& 1.7012 & 1.8246 & 1.8597 \\  
\hline    $\lambda_3 -\lambda_{3,h}$ &  3.1020e+1& 1.8466e+1 &7.0955e+0& 2.1515e+0 & 6.0054e-1 & 1.6559e-1  \\ 
\hline   order   &  & 0.7483 & 1.3798&  1.7215 & 1.8410 & 1.8586  \\  
\hline    $\lambda_4 -\lambda_{4,h}$  & 3.7904e+1& 2.3256e+1 & 9.3235e+0& 2.8817e+0 & 8.0226e+-1 & 2.1629e-1\\ 
\hline   order  &  &0.7047 &1.3186& 1.6939 & 1.8447 & 1.8910  \\  
\hline    $\lambda_5 -\lambda_{5,h}$  & 4.4211e+1& 2.8086e+1 & 1.1827e+1& 3.8023e+0 & 1.1078e+0 & 3.2207e-1  \\ 
\hline   order  &  & 0.6545 &1.2477& 1.6371 & 1.7791 &1.7823 \\

\hline
\end{tabular}   
}
\end{center}   
\end{table}

\begin{figure}[H]
\centering
\subfigure[$\gamma(h)=C$]{\includegraphics[scale=0.4]{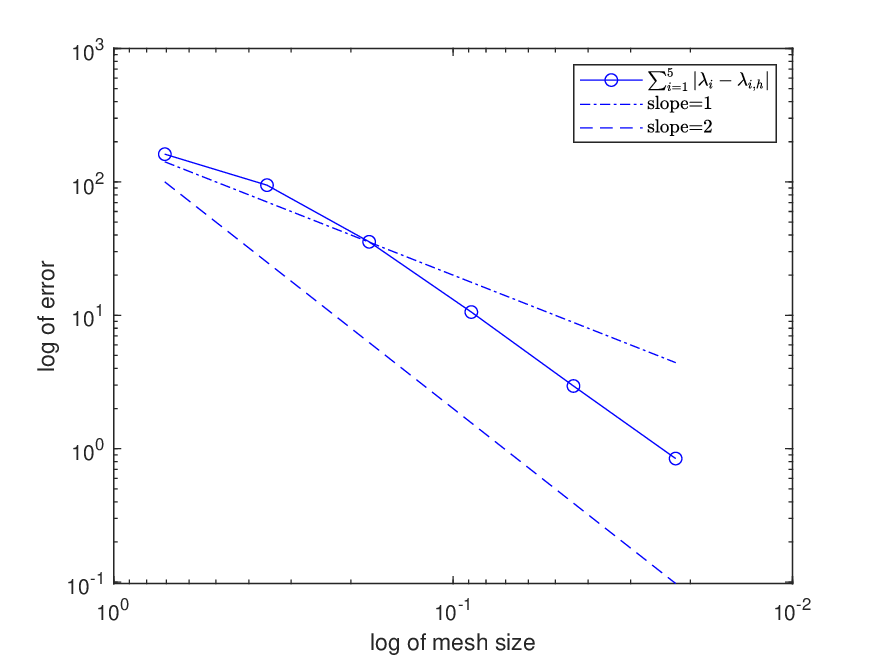}}
\subfigure[$\gamma(h)=h^{0.1}$]{\includegraphics[scale=0.4]{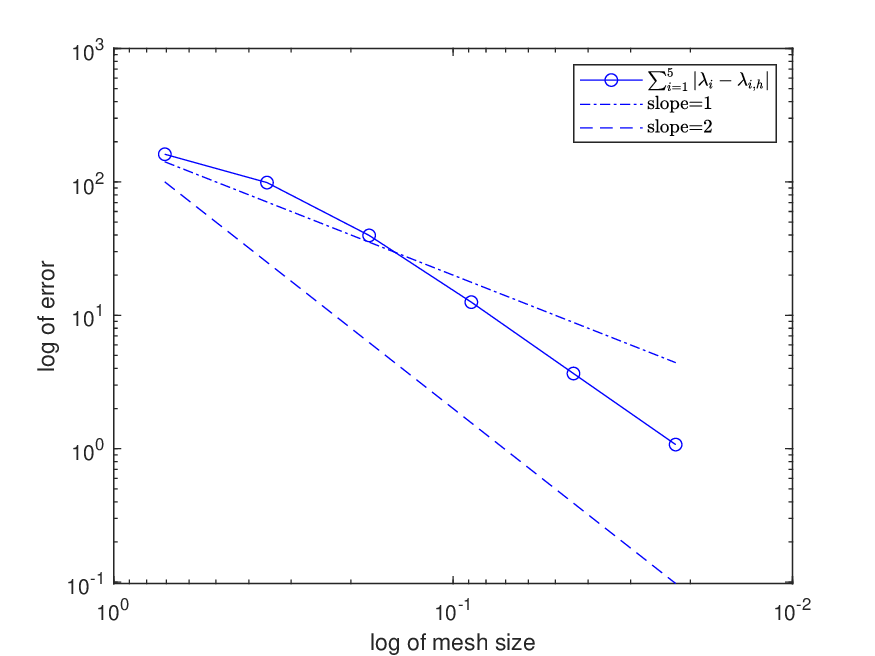}}
\subfigure[$\gamma(h)=-\frac{1}{log(h)}$]{\includegraphics[scale=0.4]{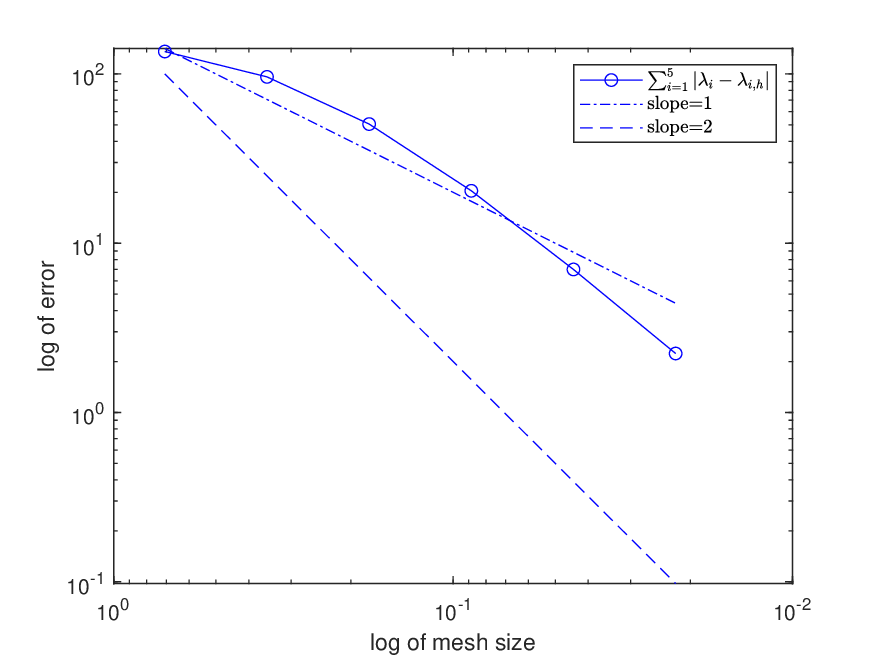}}
\caption{error estimate of eigenvalue in L-shaped region }
\end{figure}

\begin{figure}[H]
\centering
\subfigure[]{\includegraphics[scale=0.40]{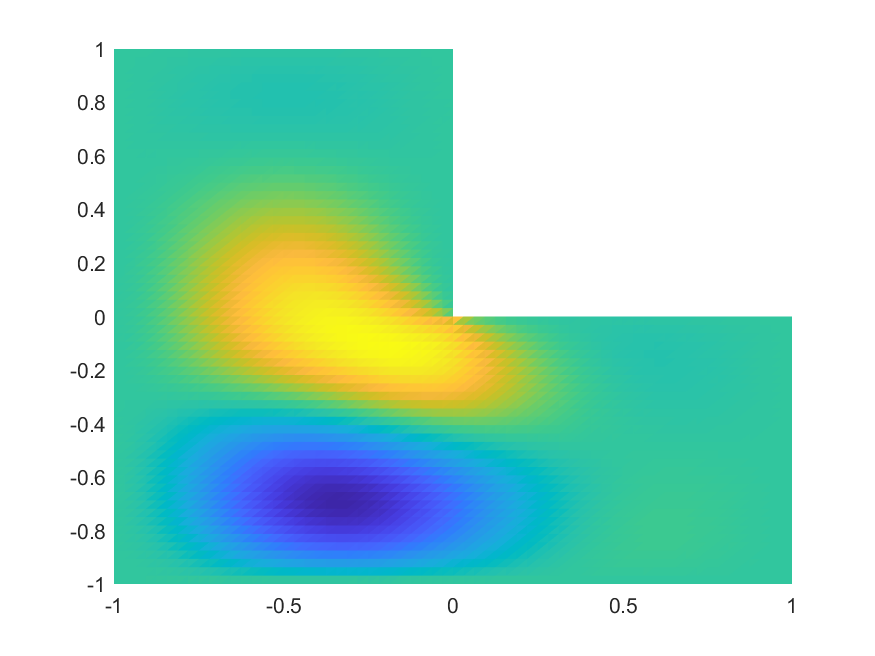}}
\subfigure[]{\includegraphics[scale=0.40]{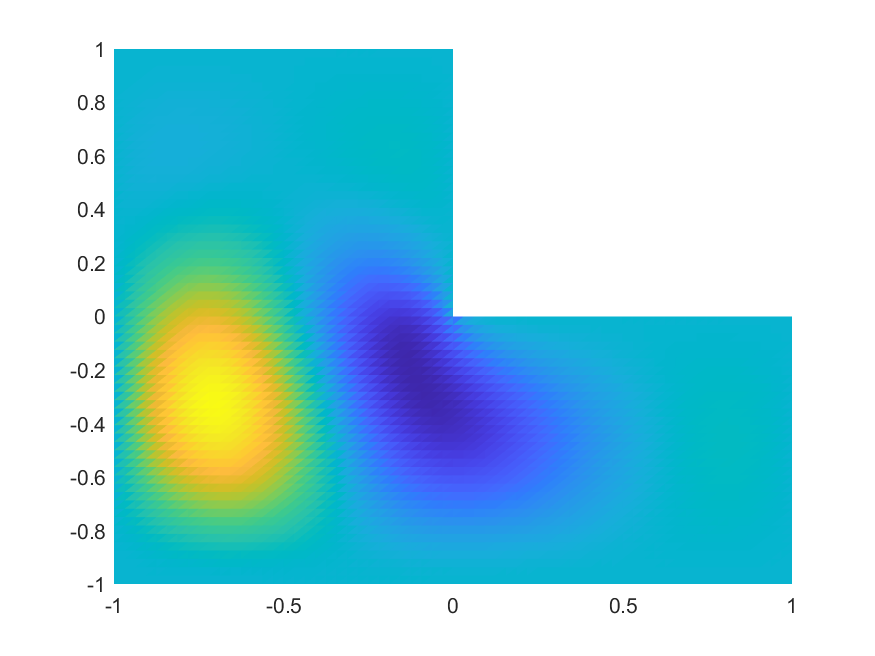}}
\caption{$(a)\text{eigenfunction }\mathbf{u}_1,(b)\text{eigenfunction }\mathbf{u}_2$}
\end{figure}

\begin{figure}[H]
\centering
\subfigure[]{\includegraphics[scale=0.42]{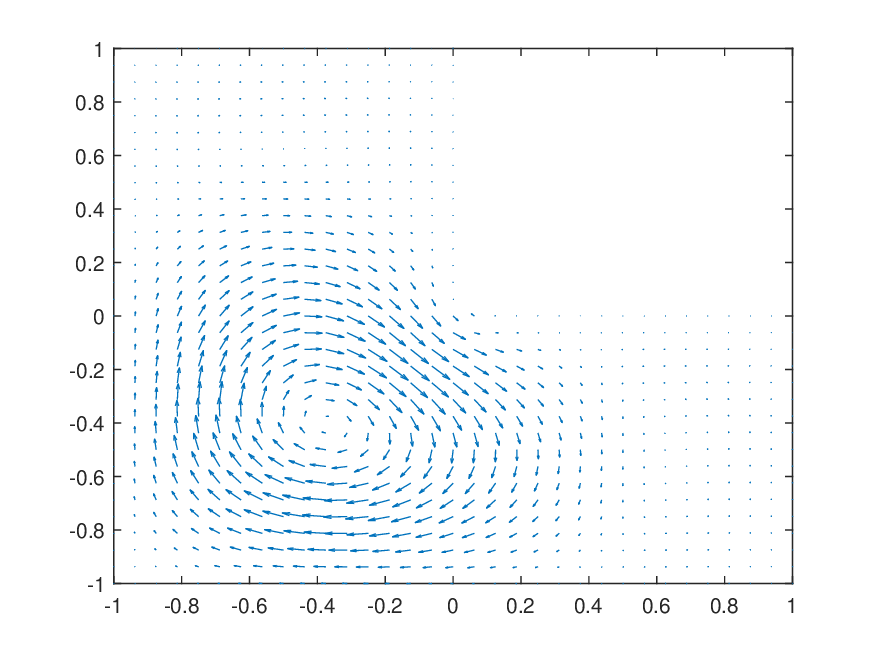}}
\subfigure[]{\includegraphics[scale=0.42]{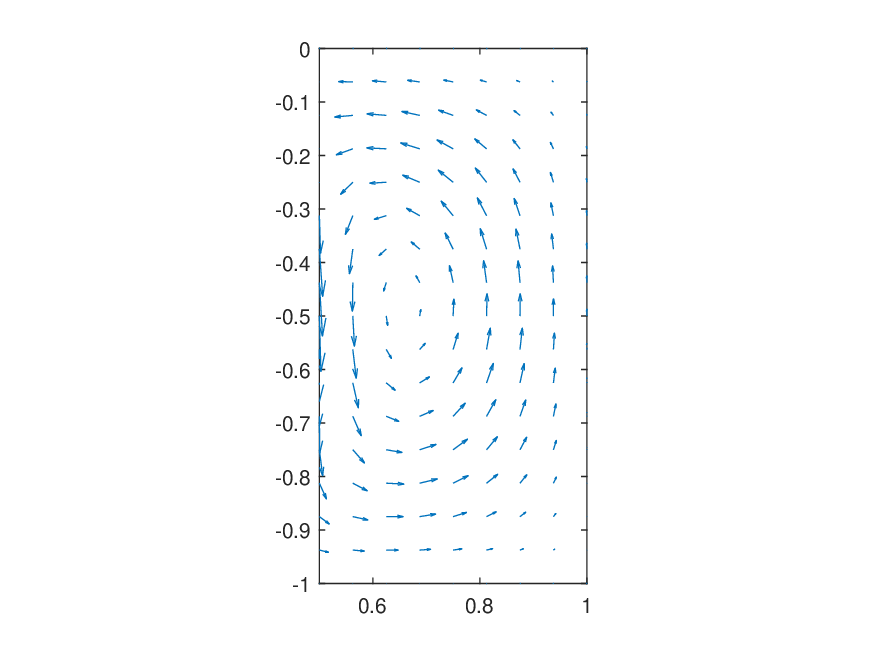}}
\subfigure[]{\includegraphics[scale=0.5]{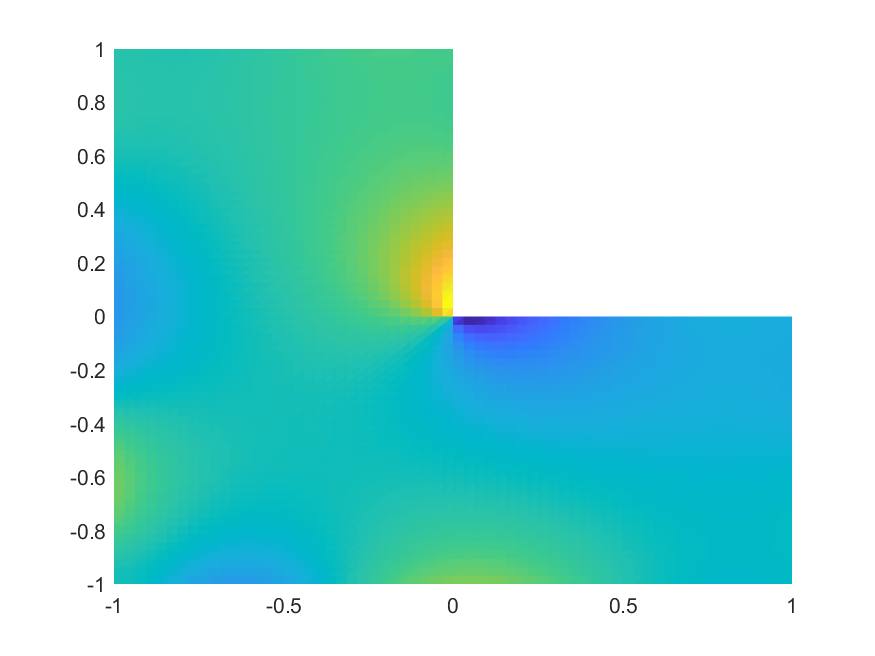}}
\caption{$(a)\text{eigenfunction }\mathbf{u} (b)\text{eigenfunction }\mathbf{u}\text{ lower right section },(c)\text{eigenfunction }p$}
\end{figure}

\section{Concluding remarks and ongoing work}
In this paper, we give the numerical scheme, error estimates, asymptotic lower bound estimates and GLB properties of the WG method for the Stokes eigenvalue problem, and verify its correctness by numerical experiments. Then WG method allows us to solve the Stokes eigenvalue problem efficiently and gives asymptotic lower bound estimates for the eigenvalues naturally without the need for a posteriori error estimates and without constructing a special finite element scheme, which is a great improvement compared to the conforming and complex non-conforming methods. At the same time, this method can be easily extended to general polygons or polyhedra in the three-dimensional region, which is more robust. In addition, the nature of WG method also allows us to perform numerical experiments on general polygonal meshes, which can obtain richer numerical results and also make this method have more application scenarios. 

There are two main difficulties in the eigenvalue problem of PDE, one is how to obtain the lower bound of the eigenvalue, and the other is the huge computational effort in solving the eigenvalue problem numerically. In this paper, the first problem is solved by the WG method, but the study of efficient numerical algorithms for the eigenvalue problem still needs to do more, in which the acceleration of the weak finite element method will be the main focus. In our future work, we will further explore and study the acceleration of the WG method using two-space, two-grid method and so on.

\section{appendix A}
In this section, we will introduce some of the technical tools used in the previous sections.
Consider Stokes equation
\begin{equation}	
    \begin{cases}
        -\Delta \mathbf{u}+\nabla p=\mathbf{f},& \text{in }\Omega,\\
        \nabla \cdot \mathbf{u}=0,& \text{in }\Omega ,\\ 
        \mathbf{u}=\mathbf{0},& \text{on }\partial\Omega, \\
     \end{cases}
\end{equation}
where $\mathbf{f}$ is given $L^2$ function. 

The corresponding WG scheme: Find $(\mathbf{u}_h,p_h) \in  V_{h,0} \times Q_h$ such that\\
\begin{equation}
    \begin{cases}
        a_w(\mathbf{u}_h,\mathbf{v}_h)-c_w(\mathbf{v}_h,p_h)=(\mathbf{f},\mathbf{v}_0), \forall \mathbf{v}_h \in V_{h,0},\\
        c_w(\mathbf{u}_h,q_h)=0,      \forall q_h \in Q_h,
     \end{cases}
\end{equation}\\
Denote $\mathbf{u}\in [H_0^1(\Omega)]^d$ and $p\in L_0^2(\Omega)$ is the solution of problem (10.1), $\mathbf{u}_h=\{\mathbf{u}_0,\mathbf{u}_b\}\in V_{h,0}$ and $p_h \in Q_h$ is the numerical solution,  denote $\mathbf{e}_h$ and $\epsilon_h$ corresponding error, i.e,
$$\mathbf{e}_h=\{\mathbf{e}_0,\mathbf{e}_b\}=\{Q_0 \mathbf{u}-\mathbf{u}_0,Q_b \mathbf{u}-\mathbf{u}_b\},\quad \epsilon_h=\mathbb{Q}_h p-p_h.$$

\setcounter{equation}{0}
\renewcommand{\theequation}{A.\arabic{equation}}

\begin{lemma}
Denote $e_h$ and $\epsilon_h$ is numerical solution of problem (10.1), then we have
\begin{equation}
    a_w(\mathbf{e}_h,\mathbf{v})-c_w(\mathbf{v},\epsilon_h)=l_{\mathbf{u}} (\mathbf{v})-\theta_p (\mathbf{v})+s(Q_h \mathbf{u},\mathbf{v}),  \quad\forall \mathbf{v} \in V_{h,0},
\end{equation}
\begin{equation}
    c_w(\mathbf{e}_h,q)=0,\quad\forall q \in Q_h,
\end{equation}
where $l_{\mathbf{u}} (\mathbf{v})=\displaystyle\sum_{T\in T_h}\left \langle \mathbf{v}_0-\mathbf{v}_b,\nabla \mathbf{u}\cdot \mathbf{n}-Q_h(\nabla \mathbf{u})\cdot \mathbf{n}\right \rangle_{\partial T},\theta_p(v)=\displaystyle\sum_{T\in T_h}\left \langle \mathbf{v}_0-\mathbf{v}_b,(\rho\right.$\\$\left.-\mathbb{Q}_h \rho)\mathbf{n}\right \rangle_{\partial T}$ both are linear form defined on $V_{h,0}$.
\end{lemma}

Proof: Since the first formula of equation (10.1) holds, from the lemma 5.2 of \cite{Wang2016AWG} we have,
$$(\nabla_w(Q_h \mathbf{u}),\nabla_w \mathbf{v})-(\nabla_w\cdot \mathbf{v},\mathbb{Q}_h p)=(\mathbf{f},\mathbf{v}_0)+l_{\mathbf{u}} (\mathbf{v})-\theta_p (\mathbf{v})$$
Add the stabilizer $s(Q_h \mathbf{u},\mathbf{v})$ to both sides of the equation, we have
$$a_w(Q_h \mathbf{u},\mathbf{v})-c_w(\mathbf{v},\mathbb{Q}_h p)=(\mathbf{f},\mathbf{v}_0)+l_{\mathbf{u}} (\mathbf{v})-\theta_p(\mathbf{v})+s(Q_h \mathbf{u},\mathbf{v}).$$
The difference with the first formula of (10.2) gives
$$a_w(\mathbf{e}_h,\mathbf{v})-c_w(\mathbf{v},\epsilon_h)=l_{\mathbf{u}} (\mathbf{v})-\theta_{\rho}(\mathbf{v})+s(Q_h \mathbf{u},\mathbf{v}).\quad\forall \mathbf{v}\in V_{h,0}$$
Besides, $0=(\nabla\cdot \mathbf{v},q)=(\nabla_w\cdot Q_h \mathbf{u},q).$ The difference with the first equation of (10.2) gives 
$$c_w(\mathbf{e}_h,q)=0,\forall q\in Q_h.$$
We can give the following estimation of $l_{\mathbf{w}}(\mathbf{v}),\theta_p (\mathbf{v})$ and $s(Q_h \mathbf{w},\mathbf{v})$.
\begin{lemma}
Denote $1\leq r\leq k$ and $\mathbf{w}\in [H^{r+1}(\Omega)]^d,\rho\in H^r(\Omega),\mathbf{v}\in V_{h,0}$, then the following estimate holds, 
$$|s(Q_h \mathbf{w},\mathbf{v})|\leq C\gamma(h)^{-1}h^r\parallel \mathbf{w}\parallel_{r+1}\parallel \mathbf{v}\parallel_V,$$
$$|l_{\mathbf{w}}(\mathbf{v})|\leq Ch^r\parallel \mathbf{w}\parallel_{r+1}\parallel \mathbf{v}\parallel_V,$$
$$|\theta_{\rho}(\mathbf{v})|\leq Ch^r\parallel \rho\parallel_{r}\parallel \mathbf{v}\parallel_V,$$
\end{lemma}

Proof:$|s(Q_h \mathbf{w},\mathbf{v})|=|\gamma(h)^{-1}\displaystyle\sum_{T\in T_h}h_T^{-1}\left \langle Q_b(Q_0 \mathbf{w})-Q_b \mathbf{w},Q_b \mathbf{v}_0-\mathbf{v}_b\right \rangle_{\partial T}|$\\
$=|\gamma(h)^{-1}\displaystyle\sum_{T\in T_h}h_T^{-1}\left \langle Q_b(Q_0 \mathbf{w}-\mathbf{w})-Q_b \mathbf{w},Q_b \mathbf{v}_0-\mathbf{v}_b\right \rangle_{\partial T}|$\\
$=|\gamma(h)^{-1}\displaystyle\sum_{T\in T_h}h_T^{-1}\left \langle Q_0 \mathbf{w}-\mathbf{w},Q_b \mathbf{v}_0-\mathbf{v}_b\right \rangle_{\partial T}|$\\
$\leq \gamma(h)^{-1}(\displaystyle\sum_{T\in T_h}(h_T^{-2}\parallel Q_0 \mathbf{w}-\mathbf{w}\parallel_T^2+\parallel \nabla(Q_0 \mathbf{w}-\mathbf{w})\parallel_T^2))^{\frac{1}{2}}(\displaystyle\sum_{T\in T_h}h_T^{-1}\parallel Q_b \mathbf{v}_0-\mathbf{v}_b\parallel_T^2)^{\frac{1}{2}}$\\
$\leq C\gamma(h)^{-1}h^r\parallel \mathbf{w}\parallel_{r+1}\parallel \mathbf{v}\parallel_V$\\
$|l_{\mathbf{w}}(\mathbf{v})|=|\displaystyle\sum_{T\in T_h}\left \langle \mathbf{v}_0-\mathbf{v}_b,\nabla \mathbf{w}\cdot \mathbf{n}-Q_h(\nabla \mathbf{w})\cdot \mathbf{n}\right \rangle_{\partial T}|\leq|\!\displaystyle\sum_{\mathcal{T\in T_h}}\left \langle \mathbf{v}_0-Q_b \mathbf{v}_0,\nabla \mathbf{w}\cdot \mathbf{n}\right.$\\$\left.-Q_h(\nabla \mathbf{w})\cdot \mathbf{n}\right \rangle_{\partial T}|+|\!\displaystyle\sum_{\mathcal{T\in T_h}}\left \langle Q_b \mathbf{v}_0-\mathbf{v}_b,\!\nabla \mathbf{w}\cdot \mathbf{n}-Q_h(\nabla \mathbf{w})\cdot \mathbf{n}\right \rangle_{\partial T}|$\\

We estimate those two parts respectively,\\
$|\displaystyle\sum_{T\in T_h}\left \langle \mathbf{v}_0-Q_b \mathbf{v}_0,\nabla \mathbf{w}\cdot \mathbf{n}-Q_h(\nabla \mathbf{w})\cdot \mathbf{n}\right \rangle_{\partial T}|$\\
$\leq C\displaystyle\sum_{T\in T_h}h_T\parallel\nabla \mathbf{w}\cdot \mathbf{n}-Q_h(\nabla \mathbf{w})\cdot \mathbf{n}\parallel_{\partial T}\parallel \nabla \mathbf{v}_0\parallel_{\partial T}$\\
$\leq C(\displaystyle\sum_{T\in T_h}h_T\parallel\nabla \mathbf{w}\cdot \mathbf{n}-Q_h(\nabla \mathbf{w})\cdot \mathbf{n}\parallel_{\partial T}^2 )^{\frac{1}{2}}(\displaystyle\sum_{T\in T_h}\parallel \nabla \mathbf{v}_0\parallel_{\partial T}^2)^{\frac{1}{2}}$\\
$\leq Ch^r \parallel \mathbf{v}\parallel_{r+1}\parallel \mathbf{v}\parallel_V$\\
$|\displaystyle\sum_{T\in T_h}\left \langle Q_b v_0-v_b,\nabla \mathbf{w}\cdot \mathbf{n}-Q_h(\nabla \mathbf{w})\cdot \mathbf{n}\right \rangle_{\partial T}|$\\
$\leq C(\displaystyle\sum_{T\in T_h}h_T\parallel\nabla \mathbf{w}\cdot \mathbf{n}-Q_h(\nabla \mathbf{w})\cdot \mathbf{n}\parallel_{\partial T}^2 )^{\frac{1}{2}}(\displaystyle\sum_{T\in T_h}h_T^{-1}\parallel Q_b \mathbf{v}_0-\mathbf{v}_b\parallel_{\partial T}^2)^{\frac{1}{2}}$\\
$\leq Ch^r\parallel \mathbf{w}\parallel_{r+1}\parallel \mathbf{v}\parallel_V.$\\
$|\theta_{\rho}(\mathbf{v})|=|\displaystyle\sum_{T\in T_h}\left \langle \mathbf{v}_0-\mathbf{v}_b,(\rho-\mathbb{Q}_h \rho)\mathbf{n}\right \rangle_{\partial T}|$
$\leq Ch^r\parallel\rho\parallel_r \parallel \mathbf{v}\parallel_V$\\
After completing the above error estimates, we can obtain the error estimates for WG method of the Stokes equation.
\begin{theorem}
Denote $(\mathbf{u},p)\in [H_0^1(\Omega) \cap H^{k+1}(\Omega)]^d\times (L_0^2(\Omega)\cap H^k(\Omega))$ and $(\mathbf{u}_h,p_h)\in V_{h,0}\times Q_h$ are the exact solution and the corresponding numerical solution of problem (10.1), respectively. Then the following error estimate holds:
$$\parallel Q_h \mathbf{u}-\mathbf{u}_h\parallel_V +\parallel \mathbb{Q}_h p-p_h\parallel\leq C\gamma(h)^{-1}h^k(\parallel \mathbf{u}\parallel_{k+1}+\parallel p\parallel_{k}).$$
\end{theorem}
Proof: Take $\mathbf{v}=\mathbf{e}_h$ into (A.1) and take $q=\epsilon_h$ into (A.2), add both two formula, then we have
$$\parallel \mathbf{e}_h\parallel_V=(\mathbf{f},\mathbf{e}_h)+l_{\mathbf{u}} (\mathbf{e}_h)-\theta_p(\mathbf{e}_h)+s(Q_h \mathbf{u},\mathbf{e}_h)$$
Hence,
$$\parallel \mathbf{e}_h\parallel_V\leq C\gamma(h)^{-1}h^k(\parallel \mathbf{u}\parallel_{k+1}+\parallel p\parallel_k),$$
And $b(\mathbf{v},\epsilon_h)=a(\mathbf{e}_h,v)-(\mathbf{f},\mathbf{v}_0)-l_{\mathbf{u}} (\mathbf{v})+\theta_p(\mathbf{v})-s(Q_h \mathbf{u},\mathbf{v}),$, from the above estimation, we can obtain,
$$|b(\mathbf{v},\epsilon_h)|\leq C\gamma(h)^{-1}h^k(\parallel \mathbf{u}\parallel_{k+1}+\parallel p\parallel_k)\parallel \mathbf{v}\parallel_V.$$
Combing the inf-sup condition of Stokes equation, we obtain
$$\parallel \epsilon_h\parallel\leq C\gamma(h)^{-1}h^k(\parallel \mathbf{u}\parallel_{k+1}+\parallel p\parallel_k),$$
By classic Nistche technique, we can obtain its $L^2$ estimation.
\begin{theorem}
Denote $(\mathbf{u},p)\in [H_0^1(\Omega) \cap H^{k+1}(\Omega)]^d\times ((L_0^2)(\Omega)\cap H^k(\Omega))$ and $(\mathbf{u}_h,p_h)\in V_{h,0}\times Q_h$ are the exact solution and the corresponding numerical solution of problem (10.1), respectively. Then the following error estimate holds$$\parallel Q_h \mathbf{u}-\mathbf{u}_h\parallel_V +\parallel \mathbb{Q}_h p-p_h\parallel\leq C\gamma(h)^{-1}h^{k+1}(\parallel \mathbf{u}\parallel_{k+1}+\parallel p\parallel_{k}).$$
\end{theorem}

----------------------------------------------

\bibliographystyle{siam}
\bibliography{library}

\end{document}